\documentclass[12pt,reqno,a4wide]{article}
\usepackage{eurosym}
\usepackage{a4wide}
\usepackage{comment}
\usepackage{mathptmx}
\usepackage{amssymb}
\usepackage{amsfonts}
\usepackage{amsthm}
\usepackage{amsmath}
\usepackage{dsfont}
\usepackage{graphicx}
\usepackage{shadow}
\usepackage[all]{xy}
\usepackage{enumerate}
\usepackage{makeidx}
\usepackage{mathrsfs}
\usepackage{upgreek}
\usepackage{stmaryrd}
\usepackage[utf8]{inputenc}
\usepackage{array}
\usepackage{tikz-cd}
\usepackage{thumbpdf}
\usepackage[colorlinks=true,pagebackref]{hyperref}
\usepackage{textcomp}
\usepackage{chapterbib}
\usepackage{graphics}
\usepackage{longtable}
\usepackage{epsf}
\usepackage{bm}
\usepackage{adjustbox}
\usepackage{centernot}

\makeindex
\newtheorem{theorem}{Theorem}[section]
\newtheorem{lem}[theorem]{Lemma}
\newtheorem{thm}[theorem]{Theorem}
\newtheorem{prop}[theorem]{Proposition}
\newtheorem{cor}[theorem]{Corollary}
\theoremstyle{definition}
\newtheorem*{Beweis}{Proof}
\newtheorem{defn}[theorem]{Definition}
\newtheorem{definition}[theorem]{Definition}

\newtheorem{rem}[theorem]{Remark}

\newtheorem{punto}[theorem]{}

\theoremstyle{remark}

\newtheorem{ex}[theorem]{Example}
\newtheorem{exs}[theorem]{Examples}

\input{tcilatex}
\begin{document}

\title{Injective Semimodules - Revisited\thanks{%
MSC2010: Primary 16Y60; Secondary 16D50 \newline
Key Words: Semirings; Semimodules; Injective Semimodules; Exact Sequences
\newline
The authors would like to acknowledge the support provided by the Deanship
of Scientific Research (DSR) at King Fahd University of Petroleum $\&$
Minerals (KFUPM) for funding this work through projects No. RG1304-1 $\&$
RG1304-2}}
\author{$%
\begin{array}{ccc}
\text{Jawad Abuhlail}\thanks{\text{Corresponding Author}} &  & \text{Rangga
Ganzar Noegraha}\thanks{\text{The paper is extracted from his Ph.D.
dissertation under the supervision of Prof. Jawad Abuhlail.}} \\
\text{abuhlail@kfupm.edu.sa} &  & \text{rangga.gn@universitaspertamina.ac.id}
\\
\text{Department of Mathematics and Statistics} &  & \text{Universitas
Pertamina} \\
\text{King Fahd University of Petroleum $\&$ Minerals} &  & \text{Jl. Teuku
Nyak Arief} \\
\text{31261 Dhahran, KSA} &  & \text{Jakarta 12220, Indonesia}%
\end{array}%
$}
\date{\today }
\maketitle

\begin{abstract}
Injective modules play an important role in characterizing different classes
of rings (e.g. Noetherian rings, semisimple rings). Some semirings have no
non-zero injective semimodules (e.g. the semiring of non-negative integers).
In this paper, we study some of the basic properties of the so called $e$%
\emph{-injective semimodules} introduced by the first author using a new
notion of \textit{exact sequences} of semimodules. We clarify the
relationships between the injective semimodules, the $e$-injective
semimodule, and the $i$-injective semimodules through several implications,
examples and counter examples. Moreover, we provide partial results for the
so called \emph{Embedding Problem} (of semimodules in injective semimodules).
\end{abstract}


\section*{Introduction}

\emph{Semirings} (defined, roughly, as rings not necessarily with
subtraction) generalize both rings and distributive bounded lattices.
Semirings and their {\emph{semimodules} (defined, roughly, as modules not
necessarily with subtraction) have many applications in Mathematics,
Computer Science and Theoretical Science (e.g., \cite{HW1998}, \cite{Gla2002}%
, \cite{LM2005}). Our main reference on semirings and their applications is
Golan's book \cite{Gol1999}, and Our main reference in rings and modules is
\cite{Wis1991}}.

\bigskip

The notion of injective objects can be defined in any category relative to a
suitable \emph{factorization system.} Injective semimodules have been
studied intensively (see \cite{Gla2002} for details). {Recently, }several
papers established homological characterizations of special classes of
semirings using (cf., \cite{KNT2009}, \cite{Ili2010}, \cite{KN2011}, \cite%
{Abu2014}, \cite{KNZ2014}, \cite{AIKN2015}, \cite{IKN2017}, \cite{AIKN2018}%
). For example, left (right) $V$-semirings, all of whose \emph{%
congruence-simple} left (right) semimodules are injective have been
completely characterized in \cite{AIKN2015}.

\bigskip

In addition to the \emph{classical notions} of \emph{injective semimodules}
over a semiring, several other notions were considered in the literature,
e.g. the so called $i$-injective\emph{\ semimodules} \cite{Alt2003} and the $%
k$\emph{-injective semimodules }\cite{KNT2009}. One reason for the interest
of such notions is the phenomenon that assuming that \emph{all} semimodules
of a given semiring $S$ to be injective forces the semiring to be a
(semisimple) \emph{ring} (cf. \cite[Theorem 3.4]{Ili2010}). Using a new
notion of exact sequences of semimodules over a semiring, Abuhlail \cite%
{Abu2014-CA} introduced a \emph{homological notion} of \emph{exactly
injective semimodules (}$e$\emph{-injective} \emph{semimodules} for short)
assuming that an appropriate $Hom$ functors preserve short exact sequences.
Such semimodules were called initially \emph{uniformly injective semimodules}
and used in \cite{Abu2014-SF} under the name \emph{normally injective
semimodules}; the terminology $e$-injective semimodules was used first in
\cite{AIKN2018}.

\bigskip

The paper is divided into three sections.

\bigskip

In Section 1, we collect some basic definitions, examples and preliminaries
used in this paper. In particular, we recall the definition and basic
properties of \emph{exact sequences} in the sense of Abuhlail \cite{Abu2014}.

\bigskip

In Section 2, we investigate mainly the $e$\emph{-injective semimodules}
over a semiring and clarify their relationships with the \emph{injective
semimodules} and the $i$\emph{-injective semimodules}. In Lemma \ref{ret-inj}
and Proposition \ref{prod-inj} we provide homological detailed proofs of the
fact that the class of injective left semimodules is closed under retracts
and direct products. It was shown in \cite[Proposition-Example 4.6.]%
{AIKN2018} that, for an additively idempotent division semiring $D,$ the
class of $e$-injective $D$-semimodules is \emph{strictly} larger than the
class of injective $D$-semimodules. Subsection \ref{sub-k-not-e-inj} is
devoted to showing that for the semiring $S:=M_{2}(\mathbb{R}^{+})$, the
class of $S$-$i$-injective left semimodules is strictly larger than the
class of $S$-$e$-injective left $S$-semimodules: Lemma \ref{M2R-all-i-inj}
shows that \emph{all} left $S$-semimodules are $S$-$I$-injective, while
Example \ref{k-inj-not-e-ink} provides a left $S$-semimodule which is not $S$%
-$e$-injective.

\bigskip

In Section 3, we investigate the so called \emph{Embedding Problem}. While
every module over a ring $R$ can be embedded in an injective semimodules,
and a module $M$ is injective if $M$ is $R$-injective (using the \emph{%
Baer's Criterion}), any semiring whose category of semimodules has both of
these nice properties is a \emph{ring} \cite[Theorem 3]{Ili2008}. Call a
left $S$-semimodule $c$-$i$-\emph{injective} if it is $M$-$i$-injective for
every \emph{cancellative} left $S$-semimodule $M.$ We prove in Theorem \ref%
{c-i-injective} that every left $S$-semimodule can be embedded \emph{%
subtractively} in a $c$-$i$-\emph{injective} left $S$-semimodule.

\section{Preliminaries}

\label{prelim}

\qquad In this section, we provide the basic definitions and preliminaries
used in this work. Any notions that are not defined can be found in {our
main reference \cite{Gol1999}. We refer to \cite{Wis1991} for the
foundations of Module and Ring Theory.}

\begin{definition}
(\cite{Gol1999}) A \textbf{semiring}%
\index{Semiring} is a datum $(S,+,0,\cdot ,1)$ consisting of a commutative
monoid $(S,+,0)$ and a monoid $(S,\cdot ,1)$ such that $0\neq 1$ and%
\begin{eqnarray*}
a\cdot 0 &=&0=0\cdot a%
\text{ for all }a\in S; \\
a(b+c) &=&ab+ac\text{ and }(a+b)c=ac+bc\text{ for all }a,b,c\in S.
\end{eqnarray*}%
A semiring $S$ with $(S,\cdot ,1)$ a commutative monoid is called a \textbf{%
commutative semiring. }A semiring $S$ with $a+a=a$ for all $a\in S$ is said
to be an \textbf{additively idempotent semiring. }A semiring with no
non-zero zero-divisors is called \textbf{entire. }We set%
\begin{equation}
V(S):=\{s\in S\mid s+t=0\text{ for some }t\in S\}.  \label{V(S)}
\end{equation}%
If $V(S)=\{0\}$, then we say that $S$ is \textbf{zerosumfree.}
\end{definition}

\begin{exs}
({\cite{Gol1999})}

\begin{itemize}
\item Every ring is a cancellative semiring.

\item Any \emph{distributive bounded lattice} $\mathcal{L}=(L,\vee ,1,\wedge
,0)$ is a commutative idempotent semiring and $1$ is an infinite element of $%
\mathcal{L}$.

\item The sets $(\mathbb{Z}^{+},+,0,\cdot ,1)$ (resp. $(\mathbb{Q}%
^{+},+,0,\cdot ,1),$ $(\mathbb{Q}^{+},+,0,\cdot ,1)$) of non-negative
integers (resp. non-negative rational numbers, non-negative real numbers) is
a commutative cancellative semiring which is not a ring.

\item $M_{n}(S),$ the set of all $n\times n$ matrices over a semiring $S,$
is a semiring.

\item $\mathbb{B}:=\{0,1\}$ with $1+1=1$ is an additively idempotent
semiring called the \textbf{Boolean semiring}.

\item $\mathbb{R}_{\max }:=(\mathbb{R}\cup \{-\infty \},\max ,-\infty ,+,0)$
is an additively idempotent semiring.
\end{itemize}
\end{exs}

\begin{punto}
\cite{Gol1999} Let $S$ and $T$ be semirings. The categories $_{S}\mathbf{SM}$
of \textbf{left} $S$-\textbf{semimodules} with morphisms the $S$-linear
maps, $\mathbf{SM}_{T}$ of right $S$-semimodules with morphisms the $T$%
-linear maps, and $_{S}\mathbf{SM}_{T}$ of $(S,T)$-bisemimodules with
morphisms the $S$-linear $T$-linear maps are defined as for left (right)
modules and bimodules over rings. The set of \emph{cancellative elements }of
a (bi)semimodule $M$ is%
\begin{equation*}
K^{+}(M):=\{m\in M\mid m+m_{1}=m+m_{2}\Longrightarrow m_{1}=m_{2}\text{ for
any }m_{1},m_{2}\in M\};
\end{equation*}%
and we say that $M$ is \textbf{cancellative,} if $K^{+}(M)=M.$ We write $%
L\leq _{S}M$ to indicate that $L$ is a \emph{subsemimodule} of the $S$%
-semimodule $M.$
\end{punto}

\begin{ex}
The category of $\mathbb{Z}^{+}$-semimodules is nothing but the category of
commutative monoids.
\end{ex}

\begin{ex}
(\cite[page 150, 154]{Gol1999}) Let $S$ be a semiring, $M$ be a left $S$%
-semimodule and $L\lvertneqq _{S}M.$ The \textbf{subtractive closure }of $L$
is defined as%
\begin{equation}
\overline{L}:=\{m\in M\mid \text{ }m+l=l^{\prime }\text{ for some }%
l,l^{\prime }\in L\}.  \label{L-s-closure}
\end{equation}%
One can easily check that $\overline{L}=Ker(M\overset{\pi }{\longrightarrow }%
M/L),$ where $\pi $ is the canonical projection. We say that $L$ is \textbf{%
subtractive,} if $L=\overline{L}.$ We say that $M$ is a \textbf{subtractive
semimodule, }if every $S$-subsemimodule $L\leq _{S}M$ is subtractive.%
\end{ex}

Following \cite{BHJK2001}, we use the following definitions.

\begin{definition}
Let $S$ be a semiring. A left $S$-semimodule $M$ is \textbf{ideal}-\textbf{%
simple}, if $0$ and $M$ are the only $S$-subsemimodules of $M$.
\end{definition}

\begin{punto}
\label{variety}(cf., \cite{AHS2004})\ The category $_{S}\mathbf{SM}$ of left
semimodules over a semiring $S$ is a closed under homomorphic images,
subobjects and arbitrary products (i.e. a \emph{variety} in the sense of
Universal Algebra). In particular, $_{S}\mathbf{SM}$ is \emph{complete},
i.e. has all limits (e.g., direct products, equalizers, kernels, pullbacks,
inverse limits) and \emph{cocomplete}, i.e. has all colimits (e.g., direct
coproducts, coequalizers, cokernels, pushouts, direct colimits). For the
construction of the pullbacks and the pushouts, see \cite{AN}.
\end{punto}

\begin{punto}
Let $M$ be a left $S$-semimodule. We say that $N\leq _{S}M$ is a

\textbf{retract} of $M$, if there exists a (surjective) $S$-linear map $%
\theta :M\longrightarrow N$ and an (injective) $S$-linear map $\psi
:N\longrightarrow M$ such that $\theta \circ \psi =\mathrm{id}_{N}$;

\textbf{direct summand }of $M,$ if there exists $L\leq _{S}M$ such that $%
M=L\oplus N.$
\end{punto}

\subsection*{Exact Sequences}

\bigskip

Throughout, $(S,+,0,\cdot ,1)$ is a semiring and, unless otherwise
explicitly mentioned, an $S$-module is a \emph{left }$S$-semimodule.

\bigskip

\begin{definition}
A morphism of left $S$-semimodules $f:L\rightarrow M$ is

$k$-\textbf{normal}, if whenever $f(m)=f(m^{\prime })$ for some $m,m^{\prime
}\in M,$ we have $m+k=m^{\prime }+k^{\prime }$ for some $k,k^{\prime }\in
Ker(f);$

$i$-\textbf{normal}, if $\func{Im}(f)=\overline{f(L)}$ ($:=\{m\in M|\text{ }%
m+l\in L\text{ for some }l\in L\}$).

\textbf{normal}, if $f$ is both $k$-normal and $i$-normal.
\end{definition}

\begin{rem}
Among others, Takahashi (\cite{Tak1981}) and Golan \cite{Gol1999} called $k$%
-normal (resp., $i$-normal, normal) $S$-linear maps $k$\emph{-regular}
(resp., $i$\emph{-regular}, \emph{regular}) morphisms. We changed the
terminology to avoid confusion with the regular monomorphisms and regular
epimorphisms in Category Theory which have different meanings when applied
to categories of semimodules.
\end{rem}

The following technical lemma is helpful in several proofs in this and
forthcoming related papers.

\begin{lem}
\label{i-normal}\emph{(\cite{AN})}\ Let $L\overset{f}{\rightarrow }M\overset{%
g}{\rightarrow }N$ be a sequence of semimodules.

\begin{enumerate}
\item Let $g$ be injective.

\begin{enumerate}
\item $f$ is $k$-normal if and only if $g\circ f$ is $k$-normal.

\item If $g\circ f$ is $i$-normal (normal), then $f$ is $i$-normal (normal).

\item Assume that $g$ is $i$-normal. Then $f$ is $i$-normal (normal) if and
only if $g\circ f$ is $i$-normal (normal).
\end{enumerate}

\item Let $f$ be surjective.

\begin{enumerate}
\item $g$ is $i$-normal if and only if $g\circ f$ is $i$-normal.

\item If $g\circ f$ is $k$-normal (normal), then $g$ is $k$-normal (normal).

\item Assume that $f$ is $k$-normal. Then $g$ is $k$-normal (normal) if and
only if $g\circ f$ is $k$-normal (normal).
\end{enumerate}
\end{enumerate}
\end{lem}

There are several notions of exactness for sequences of semimodules. In this
paper, we use the relatively new notion introduced by Abuhlail:

\begin{definition}
\label{Abu-exs}(\cite[2.4]{Abu2014}) A sequence
\begin{equation}
L\overset{f}{\longrightarrow }M\overset{g}{\longrightarrow }N  \label{LMN}
\end{equation}%
of left $S$-semimodules is \textbf{exact}, if $g$ is $k$-normal and $%
f(L)=Ker(g).$
\end{definition}

\begin{punto}
\label{def-exact}We call a sequence of $S$-semimodules $L\overset{f}{%
\rightarrow }M\overset{g}{\rightarrow }N$

\emph{proper-exact} if $f(L)=\mathrm{Ker}(g)$ (exact in the sense of
Patchkoria \cite{Pat2003});

\emph{semi-exact} if $\overline{f(L)}=\mathrm{Ker}(g)$ (exact in the sense
of Takahashi \cite{Tak1981});

\emph{quasi-exact} if $\overline{f(L)}=\mathrm{Ker}(g)$ and $g$ is $k$%
-normal (exact in the sense of Patil and Doere \cite{PD2006}).
\end{punto}

\begin{punto}
We call a (possibly infinite) sequence of $S$-semimodules
\begin{equation}
\cdots \rightarrow M_{i-1}\overset{f_{i-1}}{\rightarrow }M_{i}\overset{f_{i}}%
{\rightarrow }M_{i+1}\overset{f_{i+1}}{\rightarrow }M_{i+2}\rightarrow \cdots
\label{chain}
\end{equation}

\emph{chain complex} if $f_{j+1}\circ f_{j}=0$ for every $j;$

\emph{exact} (resp., \emph{proper-exact}, \emph{semi-exact, quasi-exact}) if
each partial sequence with three terms $M_{j}\overset{f_{j}}{\rightarrow }%
M_{j+1}\overset{f_{j+1}}{\rightarrow }M_{j+2}$ is exact (resp.,
proper-exact, semi-exact, quasi-exact).

A \textbf{short exact sequence}%
\index{short exact sequence} (or a \textbf{Takahashi extension}%
\index{Takahashi extension} \cite{Tak1982b}) of $S$-semimodules is an exact
sequence of the form%
\begin{equation*}
0\longrightarrow L\overset{f}{\longrightarrow }M\overset{g}{\longrightarrow }%
N\longrightarrow 0
\end{equation*}
\end{punto}

\begin{rem}
In the sequence (\ref{LMN}), the inclusion $f(L)\subseteq Ker(g)$ forces $%
f(L)\subseteq
\overline{f(L)}\subseteq Ker(g),$ whence the assumption $f(L)=Ker(g)$
guarantees that $f(L)=\overline{f(L)},$ \emph{i.e.} $f$ is $i$-normal. So,
the definition puts conditions on $f$ and $g$ that are dual to each other
(in some sense).
\end{rem}

The follows examples show some of the advantages of the new definition of
exact sequences over the old ones:

\begin{lem}
\label{exact}Let $L,M$ and $N$ be $S$-semimodules.

\begin{enumerate}
\item $0\longrightarrow L\overset{f}{\longrightarrow }M$ is exact if and
only if $f$ is injective.

\item $M\overset{g}{\longrightarrow }N\longrightarrow 0$ is exact if and
only if $g$ is surjective.

\item $0\longrightarrow L\overset{f}{\longrightarrow }M\overset{g}{%
\longrightarrow }N$ is semi-exact and $f$ is normal if and only if $L\simeq
\mathrm{Ker}(g).$

\item $0\longrightarrow L\overset{f}{\longrightarrow }M\overset{g}{%
\longrightarrow }N$ is exact if and only if $L\simeq \mathrm{Ker}(g)$ and $g$
is $k$-normal.

\item $L\overset{f}{\longrightarrow }M\overset{g}{\longrightarrow }%
N\longrightarrow 0$ is semi-exact and $g$ is normal if and only if $N\simeq
M/f(L).$

\item $L\overset{f}{\longrightarrow }M\overset{g}{\longrightarrow }%
N\longrightarrow 0$ is exact if and only if $N\simeq M/f(L)$ and $f$ is $i$%
-normal.

\item $0\longrightarrow L\overset{f}{\longrightarrow }M\overset{g}{%
\longrightarrow }N\longrightarrow 0$ is exact if and only if $L\simeq
\mathrm{Ker}(g)$ and $N\simeq M/L.$
\end{enumerate}
\end{lem}

\begin{cor}
\label{M/L}The following assertions are equivalent:

\begin{enumerate}
\item $0\rightarrow L\overset{f}{\rightarrow }M\overset{g}{\rightarrow }%
N\rightarrow 0$ is an exact sequence of $S$-semimodules;

\item $L\simeq \mathrm{Ker}(g)$ and $N\simeq M/f(L)$;

\item $f$ is injective, $f(L)=\mathrm{Ker}(g),$ $g$ is surjective and ($k$%
-)normal.

In this case, $f$ and $g$ are normal morphisms.
\end{enumerate}
\end{cor}

\begin{rem}
An $S$-linear map is a monomorphism if and only if it is injective. Every
surjective $S$-linear map is an epimorphism. The converse is not true in
general; for example the embedding $\iota :\mathbb{Z}^{+}\hookrightarrow
\mathbb{Z}$ is an epimorphism of $\mathbb{Z}^{+}$-semimodules.
\end{rem}

\begin{prop}
\label{adj-lim}\emph{(}cf., \emph{\cite[Proposition 3.2.2]{Bor1994})} Let $%
\mathfrak{C},\mathfrak{D}$ be arbitrary categories and $\mathfrak{C}\overset{%
\mathcal{F}}{\longrightarrow }\mathfrak{D}\overset{\mathcal{G}}{%
\longrightarrow }\mathfrak{C}$ be functors such that $(\mathcal{F},\mathcal{G%
})$ is an adjoint pair.

\begin{enumerate}
\item $\mathcal{F}$ preserves all colimits which turn out to exist in $%
\mathfrak{C}.$

\item $\mathcal{G}$ preserves all limits which turn out to exist in $%
\mathfrak{D}.$
\end{enumerate}
\end{prop}

\begin{cor}
\label{ad-l-cor}Let $S,$ $T$ be semirings and $_{T}F_{S}$ a $(T,S)$%
-bisemimodule.

\begin{enumerate}
\item $\mathrm{Hom}_{T}(-,G):$ $_{T}\mathbf{SM}\longrightarrow $ $\mathbf{SM}%
_{S}$ converts all colimits in to limits.

\item For every family of left $T$-semimodules $\{Y_{\lambda }\}_{\Lambda },$
we have a canonical isomorphism of right $S$-semimodules%
\begin{equation*}
\mathrm{Hom}_{T}(\bigoplus\limits_{\lambda \in \Lambda }Y_{\lambda
},G)\simeq \prod\limits_{\lambda \in \Lambda }\mathrm{Hom}_{T}(Y_{\lambda
},G).
\end{equation*}

\item For any directed system of left $T$-semimodules $(X_{j},\{f_{jj^{%
\prime }}\})_{J},$ we have an isomorphism of right $S$-semimodules%
\begin{equation*}
\mathrm{Hom}_{T}(\lim_{\longrightarrow }X_{j},G)\simeq \text{ }%
\lim_{\longleftarrow }\mathrm{Hom}_{T}(X_{j},G).
\end{equation*}

\item $\mathrm{Hom}_{T}(-,G)$ converts coequalizers into equalizers;

\item $\mathrm{Hom}_{T}(-,G)$ converts cokernels into kernels.
\end{enumerate}
\end{cor}

\begin{Beweis}
By \cite{KN2011}, $(G\otimes _{S}-,\mathrm{Hom}_{T}(G,-))$ is an adjoint
pair of covariant functors, where%
\begin{equation*}
G\otimes _{S}-:\text{ }_{S}\mathbf{SM}\longrightarrow \text{ }_{T}\mathbf{SM}%
\text{ and }\mathrm{Hom}_{T}(G,-):\text{ }_{T}\mathbf{SM}\longrightarrow
\mathbf{SM}_{S}.
\end{equation*}
It follows directly from Proposition \ref{adj-lim} that $\mathcal{G}:=%
\mathrm{Hom}_{T}(G,-)$ \emph{preserves} limits, whence the \emph{%
contravariant} functor $\mathrm{Hom}_{T}(-,G):$ $_{T}\mathbf{SM}%
\longrightarrow $ $\mathbf{SM}_{S}$ \emph{converts} colimits to limits. In
particular, $\mathrm{Hom}_{T}(-,G)$ converts direct coproducts (resp.
coequalizers, cokernels, pushouts, direct colimits) to direct products
(resp. equalizers, kernels, pullbacks, inverse limits).$\blacksquare $
\end{Beweis}

Corollary \ref{ad-l-cor} allows us to improve \cite[Theorem 2.6]{Tak1982a}.

\begin{prop}
\label{ll-exact}Let $_{T}G_{S}$ be a $(T,S)$-bisemimodule and consider the
functor $\mathrm{Hom}_{T}(-,G):$ $_{T}\mathbf{SM}\longrightarrow \mathbf{SM}%
_{S}.$ Let%
\begin{equation}
L\overset{f}{\rightarrow }M\overset{g}{\rightarrow }N\longrightarrow 0
\label{L-l}
\end{equation}%
be a sequence of left $T$-semimodules and consider the sequence of right $S$%
-semimodules%
\begin{equation}
0\longrightarrow \mathrm{Hom}_{T}(N,G)\overset{(g,G)}{\rightarrow }\mathrm{%
Hom}_{T}(M,G)\overset{(f,G)}{\longrightarrow }\mathrm{Hom}_{T}(L,G).
\label{LG}
\end{equation}

\begin{enumerate}
\item If $M\overset{g}{\rightarrow }N\longrightarrow 0$ is exact and $g$ is
normal, then $0\longrightarrow \mathrm{Hom}_{T}(N,G)\overset{(g,G)}{%
\rightarrow }\mathrm{Hom}_{T}(M,G)$ is exact and $(g,G)$ is normal.

\item If \emph{(\ref{L-l}) }is semi-exact and $g$ is normal, then \emph{(\ref%
{LG}) }is proper-exact \emph{(}semi-exact\emph{)} and $(g,G)$ is normal.

\item If \emph{(\ref{L-l}) }is exact and $\mathrm{Hom}_{T}(-,G)$ converts $i$%
-normal morphisms into $k$-normal ones, then \emph{(\ref{LG}) }is exact.
\end{enumerate}
\end{prop}

\begin{Beweis}
\begin{enumerate}
\item The following implications are clear: $M\overset{g}{\rightarrow }%
N\longrightarrow 0$ is exact $\Longrightarrow $ $g$ is surjective $%
\Longrightarrow $ $(g,G)$ is injective $\Longrightarrow $ $0\longrightarrow
\mathrm{Hom}_{T}(N,G)\overset{(g,G)}{\rightarrow }\mathrm{Hom}_{T}(M,G)$ is
exact. Assume that $g$ is normal and consider the exact sequence of $S$%
-semimodules%
\begin{equation*}
0\longrightarrow \mathrm{Ker}(g)\overset{\iota }{\longrightarrow }M\overset{g%
}{\longrightarrow }N\longrightarrow 0.
\end{equation*}%
Notice that $N\simeq \mathrm{Coker}(\iota ).$ By Corollary \ref{ad-l-cor}, $%
\mathrm{Hom}_{T}(-,G)$ converts cokernels into kernels, we conclude that $%
(g,G)=\mathrm{ker}((f,G))$ whence normal.

\item Apply Lemma \ref{exact} (5): $L\overset{f}{\rightarrow }M\overset{g}{%
\rightarrow }N\longrightarrow 0$ is semi-exact and $g$ is normal $%
\Longleftrightarrow $ $M\simeq \mathrm{Coker}(f).$ Since the contravariant
functor $\mathrm{Hom}_{T}(-,G)$ converts cokernels into kernels, it follows
that $\mathrm{Hom}_{T}(N,G)=\mathrm{Ker}((f,G))$ which is in turn equivalent
to (\emph{\ref{LG}}) being semi-exact and $(g,G)$ being normal. Notice that
\begin{equation*}
(g,G)(\mathrm{Hom}_{S}(N,G))=\overline{(g,G)(\mathrm{Hom}_{S}(N,G))}=\mathrm{%
Ker}((f,G)),
\end{equation*}%
\emph{i.e.} (\emph{\ref{LG}}) is proper-exact (whence semi-exact).

\item This follows immediately from \textquotedblleft 2\textquotedblright\
and the assumption on $\mathrm{Hom}_{T}(-,G).\blacksquare $
\end{enumerate}
\end{Beweis}

\section{Injective Semimodules}

\markright{\scriptsize\tt Chapter \ref{sec-inj}: $e$-injective
Semimodules} 

There are several notions of injectivity for a semimodules $M$ over a
semiring $S$ which coincide if it were a module over a ring. In this
section, we consider some of these and clarify the relationships between
them. In particular, we investigate the so called $e$\emph{-injective
semimodules} which turn to coincide with the so called \emph{normally
injective semimodules} (both notions introduced by Abuhlail and called \emph{%
uniformly injective semimodules} in \cite[1.25, 1.24]{Abu2014-CA}, the
terminology \textquotedblleft $e$\emph{-injective}" was first used in \cite%
{AIKN2018}. We also clarify their relations ships with \emph{injective
semimodules} \cite{Gol1999} and $i$\emph{-injective semimodules} \cite%
{Alt2003}.

\bigskip

As before, $(S,+,0,\cdot ,1)$ is a semiring and, unless otherwise explicitly
mentioned, and $S$-module is a left $S$-semimodule. Exact sequences here are
in the sense of Abuhlail \cite{Abu2014} (see Definition \ref{Abu-exs}).

\begin{definition}
(\cite[1.24]{Abu2014-CA}) Let $M$ be a left $S$-semimodule. A left $S$%
-semimodules $J$ is $M$-$e$-\textbf{injective}, if the contravariant functor%
\begin{equation*}
Hom_{S}(-,J):\text{ }_{S}\mathbf{SM}\longrightarrow \text{ }_{\mathbb{Z}^{+}}%
\mathbf{SM}
\end{equation*}%
transfers every short exact sequence of left $S$-semimodules
\begin{equation*}
0\longrightarrow L\overset{f}{\longrightarrow }M\overset{g}{\longrightarrow }%
N\longrightarrow 0
\end{equation*}%
into a short exact sequence of commutative monoids%
\begin{equation*}
0\longrightarrow Hom_{S}(N,J)\longrightarrow Hom_{S}(M,J)\longrightarrow
Hom_{S}(L,J)\longrightarrow 0.
\end{equation*}%
We say that $J$ is $e$\textbf{-injective,} if $J$ is $M$-$e$-injective for
every left $S$-semimodule $M.$
\end{definition}

\begin{punto}
Let $I$ be a left $S$-semimodule.

For a left $S$-semimodule $M,$ we say that $I$ is

$M$\textbf{-injective}
\index{Semimodule!injective}\cite[page 197]{Gol1999} if for every \emph{%
injective} $S$-linear map $f:L\rightarrow M$ and any $S$-linear map $%
g:L\rightarrow I,$ there exists an $S$-linear map $h:M\rightarrow I$ such
that $h \circ f=g;$%
\begin{equation*}
\xymatrix{ 0\ar[r] & L \ar[r]^{f} \ar[d]_{g} & M \ar@{.>}[ld]^{h} \\ & I }
\end{equation*}

$M$-$i$\textbf{-injective}%
\index{Semimodule!$i$-injective} \cite{Alt2003} if for every \emph{normal
monomorphism} $f:L\rightarrow M$ and any $S$-linear map $g:L\rightarrow I,$
there exists an $S$-linear map $h:M\rightarrow I$ such that $h\circ f=g;$%
\begin{equation*}
\xymatrix{ 0\ar[r] & L \ar[rrr]^{f \, (normal)} \ar[d]_{g} & & & M
\ar@{.>}[llld]^{h} \\ & I }
\end{equation*}

\textbf{normally }$M$-\textbf{injective }\cite[1.24]{Abu2014-CA} if for
every \emph{normal monomorphism} $f:L\longrightarrow M$ and any $S$-linear
map $g:L\longrightarrow I,$ there exists an $S$-linear map $%
h:M\longrightarrow I$ such that $h\circ f=g$%
\begin{equation*}
\xymatrix{0 \ar[r] & L \ar[dd]_{g} \ar[rrr]^{f \, (normal)} & & & M
\ar@<-.5ex>[r]_{h_1} \ar@<.5ex>[r]^{h_2} \ar@{.>}[llldd]^{h}
\ar@/^-0.5pc/[llldd]_{h'} & I\\ \\ & I}
\end{equation*}%
and whenever an $S$-linear map $h^{\prime }:M\rightarrow I$ satisfies $%
h^{\prime }\circ f=g$, there exist $S$-linear maps $h_{1},h_{2}:M\rightarrow
I$ such that $h_{1}\circ f=0=h_{2}\circ f$ and $h+h_{1}=h^{\prime }+h_{2}$.

We say that $I$ is injective (resp., $i$-injective, normally injective) if $%
I $ is $M$-injective (resp., $M$-$i$-injective, normally $M$-projective) for
every left $S$-semimodule $M.$
\end{punto}

\begin{prop}
\label{inj-e=n}Let $I$ be a left $S$-semimodule.

\begin{enumerate}
\item Let $M$ be a left $S$-semimodule. Then $I$ is $M$-$e$-injective if and
only if $I$ is normally $M$-injective.

\item $_{S}I$ is $e$-injective if and only if $_{S}I$ is normally injective.
\end{enumerate}
\end{prop}

\begin{Beweis}
We only need to prove (1). Let $M$ be a left $S$-semimodule.

($\Longrightarrow $) Assume that $I$ is $M$-$e$-injective. Let $L\leq _{S}M$
be a subtractive $S$-subsemimodule. By Lemma \ref{exact}, we have a short
exact sequence of left $S$-semimodules%
\begin{equation}
0\longrightarrow L\overset{\iota }{\longrightarrow }M\overset{\pi }{%
\longrightarrow }M/L\longrightarrow 0  \label{LM-inj}
\end{equation}%
where $\iota $ is the canonical embedding and $\pi $ is the canonical
projection. By our assumption, the contravariant functor $Hom_{S}(-,I):$ $%
_{S}\mathbf{SM}\longrightarrow $ $_{\mathbb{Z}^{+}}\mathbf{SM}$ preserves
exact sequences, whence the following sequence of commutative monoids%
\begin{equation*}
0\longrightarrow Hom_{S}(M/L,I)\overset{(\pi ,I)}{\longrightarrow }%
Hom_{S}(M,I)\overset{(\iota ,I)}{\longrightarrow }Hom_{S}(L,I)%
\longrightarrow 0
\end{equation*}%
is exact. In particular, $(\iota ,I):Hom_{S}(M,I)\longrightarrow
Hom_{S}(L,I) $ is a normal epimorphism, \emph{i.e.} $I$ is normally $M$%
-injective.

($\Longleftarrow $) Let%
\begin{equation}
0\longrightarrow L\overset{f}{\longrightarrow }M\overset{g}{\longrightarrow }%
N\longrightarrow 0  \label{lmn}
\end{equation}%
be an exact sequence of left $S$-semimodules. Applying the contravariant
functor $Hom_{S}(-,I)$ to (\ref{lmn}) it follows by Lemma \ref{ll-exact} (2)
and our assumption that the following sequence of commutative monoids%
\begin{equation}
0\longrightarrow Hom_{S}(N,I)\overset{(g,I)}{\longrightarrow }Hom_{S}(M,I)%
\overset{(f,I)}{\longrightarrow }Hom_{S}(L,I)\longrightarrow 0
\label{0hom-lmn}
\end{equation}%
is exact, i.e. $_{S}I$ is injective.$\blacksquare $
\end{Beweis}

The proof of the following result is similar to that of \cite[Theorem 3.7]%
{Alt2003}

\begin{prop}
\label{inj->}Let%
\begin{equation}
L\overset{f}{\longrightarrow }M\overset{g}{\longrightarrow }N  \label{lmn-33}
\end{equation}%
be a sequence of left $S$-semimodules, $I$ a left $S$-semimodule and
consider the sequence%
\begin{equation}
Hom_{S}(N,I)\overset{(g,I)}{\longrightarrow }Hom_{S}(M,I)\overset{(f,I)}{%
\longrightarrow }Hom_{S}(L,I)  \label{Hom-33}
\end{equation}
of commutative monoids.

\begin{enumerate}
\item If (\ref{lmn-33}) is exact with $g$ normal and $I$ is $i$-injective,
then (\ref{Hom-33}) is proper-exact.

\item If (\ref{lmn-33}) is exact with $g$ normal and $I$ is $e$-injective,
then (\ref{Hom-33}) is exact and $(g,I)$ is normal.

\item If (\ref{lmn-33}) is exact and $I$ is injective, then (\ref{Hom-33})
is proper exact.
\end{enumerate}
\end{prop}

\begin{Beweis}
By Corollary \ref{M/L}, we have a short exact sequence of left $S$%
-semimodules%
\begin{equation*}
0\longrightarrow Ker(g)\overset{\iota }{\longrightarrow }M\overset{\pi }{%
\longrightarrow }M/Ker(g)\longrightarrow 0
\end{equation*}%
where $\iota $ and $\pi $ are the canonical $S$-linear maps. Since (\ref%
{lmn-33}) is proper exact, $f(M)=Ker(g)$ and $M/Ker(g)=M/f(M)\simeq
Coker(f). $ By the \emph{Universal Property of Kernels}, there exists a
unique $S$-linear map $%
\widetilde{f}:L\longrightarrow Ker(g)$ such that $\iota \circ \widetilde{f}%
=f $ (and $\widetilde{f}$ is surjective). On the other hand, by the \emph{%
Universal Property of Cokernels, }there exists a unique $S$-linear map $%
\widetilde{g}:M/Ker(g)\longrightarrow N$ such that $\widetilde{g}\circ \pi
=g $. So, we have a commutative diagram of left $S$-semimodules%
\begin{equation}
\xymatrix{ & & & L \ar[lld]_{\widetilde{f}} \ar[d]^{f} & & & 0 \ar[ld] \\ 0
\ar[r] & Ker(g) \ar[ld] \ar[rr]^{\iota} & & M \ar[rr]^{\pi} \ar[d]^{g} & &
M/Ker(g) \ar[lld]^{\widetilde{g}} \ar[r] & 0 \\ 0 & & & N & & & }
\label{M-ker}
\end{equation}%
Applying the contravariant functor $Hom_{S}(-,I),$ we get the sequence%
\begin{equation}
0\longrightarrow Hom_{S}(M/Ker(g),I)\overset{(\pi ,I)}{\longrightarrow }%
Hom_{S}(M,I)\overset{(\iota ,I)}{\longrightarrow }Hom_{S}(Ker(g),I)%
\longrightarrow 0  \label{Hom-Ker}
\end{equation}%
and we obtain the commutative diagram%
\begin{equation}
\xymatrix{ & & & Hom_S(N,I) \ar[lldd]_{({\widetilde{g}},I)} \ar[dd]^{(g,I)}
& & & 0 \ar[ldd] \\ & & & & & & & \\ 0 \ar[r] & Hom_{S}(M/Ker(g),I)
\ar[rr]^{(\pi ,I)} \ar[ldd] & & Hom_{S}(M,I) \ar[rr]^{(\iota ,I)}
\ar[dd]^{(f,I)} & & Hom_{S}(Ker(g),I) \ar[lldd]^{{({\widetilde{f}},I)}}
\ar[r] & 0 \\ & & & & & & & \\ 0 & & & Hom_S(L,I) & & }  \label{homM-ker}
\end{equation}%
of commutative monoids.

\begin{enumerate}
\item Since $g$ is $k$-normal, we conclude that $\widetilde{g}$ is
injective. Moreover, $\pi $ is surjective, and $g=\widetilde{g}\circ \pi $
is normal, whence $\widetilde{g}$ is normal by Lemma \ref{i-normal} (2).
Since $_{S}I$ is $i$-injective, the sequence (\ref{Hom-Ker}) is proper-exact
and $(\widetilde{g},I)$ is surjective (see Proposition \ref{ll-exact} (2)).
It follows that%
\begin{equation*}
\begin{array}{ccccc}
Ker((f,I)) & = & Ker((\widetilde{f},I)\circ (\iota ,I)) &  &  \\
& = & Ker((\iota ,I)) &  & \text{(}(\widetilde{f},I)\text{ is injective)} \\
& = & im((\pi ,I)) &  & \text{((\ref{Hom-Ker}) is proper exact)} \\
& = & im((\pi ,I)\circ (\widetilde{g},I)) &  & \text{(}(\widetilde{g},I)%
\text{ is surjective)} \\
& = & im((g,I)). &  &
\end{array}%
\end{equation*}

\item By (1), the sequence (\ref{Hom-33})\ is proper-exact. Since $(%
\widetilde{f},I)$ is injective and $(\iota ,I)$ is $k$-normal, it follows by
Lemma \ref{i-normal} (1-a) that $(f,I)=(\widetilde{f},I)\circ (\iota ,I)$ is
$k$-normal. Consequently, (\ref{Hom-33}) is exact.

Notice that, moreover, $(\pi ,I)\ $is a normal monomorphism and $(\widetilde{%
g},I)\ $is $i$-normal, whence $(g,I)=(\pi ,I)\circ (\widetilde{g},I)$ is
normal by Lemma \ref{i-normal} (1-c).

\item The proof is similar to that of (1). Notice that by our assumption (%
\ref{lmn-33}) is exact; in particular, $g$ is $k$-normal, which is need to
show that $\widetilde{g}$ is injective, whence $(\widetilde{g},I)$ is
surjective since $_{S}I$ is injective.$\blacksquare $
\end{enumerate}
\end{Beweis}

\begin{thm}
\label{e-inj-3}Let $M$ be a left $S$-Semimodule. The following are
equivalent for a left $S$-semimodule $I:$

\begin{enumerate}
\item $_{S}I$ is normally $M$-injective;

\item $_{S}I$ is $M$-$e$-injective;

\item For every \emph{exact }sequence of left $S$-semimodules (\ref{lmn-33})
with $g$ normal, the induced sequence of commutative monoids (\ref{Hom-33})
is exact and $(g,I)$ is normal.
\end{enumerate}
\end{thm}

\begin{Beweis}
$(1)\Rightarrow (2)$ This follows by Proposition \ref{inj-e=n}.

$(2)\Rightarrow (3)$ This follows by Proposition \ref{inj->}.

$(3)\Rightarrow (1)$ This follows directly by applying the assumption to the
exact sequences of left $S$-semimodules the form $0\longrightarrow M\overset{%
g}{\longrightarrow }N$ with $g$ normal.$\blacksquare $
\end{Beweis}

Using Propositions \ref{ll-exact} and \ref{inj->}, one can easily recover
the following characterizations of injective and $i$-injective semimodules
proved in \cite{Alt2003} which inspired our characterizations of $e$%
-injective semimodules in Theorem \ref{e-inj-3}.

\begin{thm}
\label{I-inj-3}Let $M$ be a left $S$-Semimodule. The following are
equivalent for a left $S$-semimodule $I:$

\begin{enumerate}
\item $_{S}I$ is $M$-$i$-injective;

\item for every \emph{proper-exact} sequence of left $S$-semimodules (\ref%
{lmn}) in which $f$ is normal and $g$ is $k$-normal, the induced sequence of
commutative monoids (\ref{0hom-lmn}) is \emph{proper-exact};

\item for every \emph{proper-exact }sequence of left $S$-semimodules (\ref%
{lmn-33}) with $g$ normal, the induced sequence of commutative monoids (\ref%
{Hom-33}) is \emph{proper-exact}.
\end{enumerate}
\end{thm}

\begin{thm}
\label{inj-3}Let $M$ be a left $S$-Semimodule. The following are equivalent
for a left $S$-semimodule $I:$

\begin{enumerate}
\item $_{S}I$ is $M$-injective;

\item for every \emph{proper-exact} sequence of left $S$-semimodules (\ref%
{lmn}) in which $f$ is $k$-normal and $g$ is normal, the induced sequence of
commutative monoids (\ref{0hom-lmn}) is \emph{proper-exact}.

\item for every \emph{exact }sequence of left $S$-semimodules (\ref{lmn-33}%
), the induced sequence of commutative monoids (\ref{Hom-33}) is \emph{%
proper-exact}.
\end{enumerate}
\end{thm}

It follows directly from the definitions that, for any semiring $S$ and any
left $S$-semimodule $M,$ the class $\mathcal{I}_{S}^{i}(M)$ of $M$\emph{-}$i$%
\emph{-injective} left $S$-semimodules contains both the class $\mathcal{I}%
_{S}(M)$ of \emph{injective} left $S$-semimodules and the class $\mathcal{I}%
_{S}^{e}(M)$ of $S$\emph{-}$e$\emph{-injective} left $S$-semimodules, i.e.%
\begin{equation}
\mathcal{I}_{S}(M)\cup \mathcal{I}_{S}^{e}(M)\subseteq \mathcal{I}%
_{S}^{i}(M).  \label{union}
\end{equation}

While every \emph{projective} semimodule is $e$\textit{-projective (}see
\cite{AIKN2018}), it is not evident that every injective semimodule is $e$%
-injective if the base semiring is arbitrary. However, we have a partial
result:

\begin{prop}
\label{inj->e-inj}\emph{(\cite[Theorem 4.5]{AIKN2018})}\ Let $S$ be an \emph{%
additively idempotent} semiring. Then every injective left $S$-semimodule is
$e$-injective.
\end{prop}

The following examples shows that the converse of Proposition \ref%
{inj->e-inj} is not true in general.

\begin{ex}
\label{e-not-inj}(\cite[4.6]{AIKN2018}) Let $D$ be an \emph{additively
idempotent} division semiring (e.g., $D=\mathbb{B}$, the Boolean semiring).
Then $D$ has an $e$-injective left $S$-semimodule which is not injective.
\end{ex}

We illustrate first that \emph{relative injectivity} and \emph{relative }$e$%
\emph{-injectivity} are \emph{not }related even when the base semiring is
commutative and additively idempotent. The following example shows that the
relative version of Proposition \ref{inj->e-inj} is not valid, i.e. relative
injectivity does not guarantee relative $e$-injectivity.

\begin{ex}
\label{ex-inj-not-e}Consider the semiring $S:=\{0,1,a\}$ \cite{Alt2003} with
addition and multiplication given by%
\begin{equation*}
\begin{array}{cccc}
\begin{array}{cccc}
+ & 0 & 1 & a \\
0 & 0 & 1 & a \\
1 & 1 & 1 & 1 \\
a & a & 1 & a%
\end{array}
&  &  &
\begin{array}{cccc}
\cdot  & 0 & 1 & a \\
0 & 0 & 0 & 0 \\
1 & 0 & 1 & a \\
a & 0 & a & a%
\end{array}
\\
\text{addition} &  &  & \text{multiplication}%
\end{array}%
\end{equation*}%
We show that $S$ is $S$-injective but not $S$-$e$-injective. The ideals of $S
$ are $\{0\},$ $\{0,a,\}$ and $S.$ Clearly, $\{0,a\}$ is subtractive, whence
$S$ is a \emph{subtractive} semiring and our example shows that the
inclusion $\mathcal{I}_{S}^{e}(S)\subseteq \mathcal{I}_{S}^{i}(S)=\mathcal{I}%
_{S}(S)$ is strict.

\textbf{Claim I: }$S$ is $S$-injective.

We need to consider only the canonical embedding $\{0,a\}\overset{\iota }{%
\hookrightarrow }S$ and the $S$-linear map%
\begin{equation}
\varphi :\{0,a\}\longrightarrow S,\text{ }0\mapsto 0\text{ and }a\mapsto a.
\label{varphi}
\end{equation}%
Notice that $\varphi (a)\neq 1:$ if so, then $1=\varphi (a)=\varphi (a\cdot
1)=a\varphi (1)$, i.e. $a$ has a multiplicative inverse, a contradiction.
Notice that $\varphi $ can be extended to an $S$-linear map through $id_{S}.$
It follows that $S$ is $S$-injective.

\textbf{Claim II: }$S$ is \emph{not} $S$-$e$-injective.

The $S$-linear map (\ref{varphi}) can be extended through another $S$-linear
map, namely
\begin{equation*}
\widetilde{h}:S\longrightarrow S,\text{ }0\mapsto 0,\text{ }a\mapsto a,\text{
}1\mapsto a.
\end{equation*}%
However, the only $S$-linear map $h:S\longrightarrow S$ such that $(h\circ
\iota )(\{0,a\})=0$ is the $h=0:$ indeed, $h(1)=a$ implies $0=h(a)=h(a\cdot
1)=ah(1)=a\cdot a=a,$ a contradiction; and $h(1)=1$ implies $0=h(a)=h(a\cdot
1)=ah(1)=a\cdot 1=a$, a contradiction. So, we cannot find $h_{1},h_{2}\in
Hom_{S}(S,S)$ such that $id_{S}+h_{1}=\widetilde{h}+h_{2}.$ Consequently, $S$
is not $S$-$e$-injective.$\blacksquare $
\end{ex}

The following example shows that relative $e$-injectivity does not guarantee
relative injectivity. In fact, we given an example for which%
\begin{equation*}
\mathcal{I}_{S}(S)\subsetneqq \mathcal{I}_{S}^{e}(S)=\text{ }_{S}\mathbf{SM}=%
\mathcal{I}_{S}^{i}(S).
\end{equation*}

\begin{ex}
\label{ex-e-not-inj}Consider the commutative \emph{additively idempotent}
semiring $S:=(\mathbb{Z}^{+},\max ,0,\cdot ,1).$ Then $S$ has \emph{no}
non-trivial proper subtractive ideals, whence every $S$-semimodule is $S$-$e$%
-injective ($S$-$i$-injective). By \cite[Example 2.7]{Alt2003}, $S$ is not $%
S $-injective. In particular, our example shows that the inclusion $\mathcal{%
I}_{S}(S)\subseteq \mathcal{I}_{S}^{i}(S)$ is strict.
\end{ex}

Next, we provide detailed \emph{homological proofs} rather than compact
categorical ones of the facts that, for a given left $S$-semimodule $M$, the
class of $M$-$e$-injective semimodules is closed under retracts and direct
products (cf., \cite[Corollary 3.3]{AIKN2018}).

\begin{prop}
\label{ret-inj}

\begin{enumerate}
\item Let $M$ be a left $S$-semimodule. Every retract of a left $M$-$e$%
-injective $S$-semimodule is $M$-$e$-injective.

\item A retract of an $e$-injective $S$-semimodule is $e$-injective.
\end{enumerate}
\end{prop}

\begin{punto}
We need to prove (1) only.

Let $J$ be an $M$-$e$-injective left $S$-semimodule and $I$ a retract of $J$
along with $S$-linear maps $\iota :I\longrightarrow J$ and $\pi
:J\longrightarrow I$ such that $\pi \circ \iota =id_{I}.$ Let $%
f:L\rightarrow M$ be a normal $S$-monomorphism and $g:L\rightarrow I$ be an $%
S$-linear map.%
\begin{equation*}
\xymatrix{0 \ar[r] & L \ar[r]^{f} \ar[d]_{g} & M \ar@{-->}[ldd]^{h^*}\\ & I
\ar[d]_{\iota} \\ & J}
\end{equation*}%
Since $J$ is $M$-$e$-injective, there is an $S$-linear map $h^{\ast
}:M\rightarrow J$ such that $h^{\ast }\circ f=\iota \circ g$. Consider $%
h:=\pi \circ h^{\ast }.$ Then we have%
\begin{equation*}
\begin{array}{ccccc}
h\circ f & = & (\pi \circ h^{\ast })\circ f & = & \pi \circ (h^{\ast }\circ
f) \\
& = & \pi \circ (\iota \circ g) & = & (\pi \circ \iota )\circ g \\
& = & id_{I}\circ g & = & g.%
\end{array}%
\end{equation*}%
Suppose that $h^{\prime }:M\rightarrow I$ is an $S$-linear map such that $%
h^{\prime }\circ f=g.$ Notice that $\iota \circ h^{\prime }\circ f=\iota
\circ g$. Since $J$ is $M$-$e$-injective, there exist $S$-linear maps $%
h_{1}^{\ast },h_{2}^{\ast }:M\rightarrow J$ such that $h_{1}^{\ast }\circ
f=0=h_{2}^{\ast }\circ f$ and $h^{\ast }+h_{1}^{\ast }=\iota \circ h^{\prime
}+h_{2}^{\ast }$.%
\begin{equation*}
\xymatrix{0 \ar[r] & L \ar[d]_{g} \ar[r]^{f} & M \ar@<-.5ex>[r]_{h^*_1}
\ar@<.5ex>[r]^{h^*_2} \ar[ld]^{h} \ar@/^-0.5pc/[ld]_{h'}
\ar@/^0.5pc/[ldd]^{h^*} & J\\ & I \ar[d]_{\iota} \\ & J }
\end{equation*}%
Consider $h_{1}:=\pi \circ h_{1}^{\ast }$ and $h_{2}:=\pi \circ h_{2}^{\ast
}.$ Then we have, for $i=1,2$, $h_{i}\circ f=\pi \circ h_{i}^{\ast }\circ
f=\pi \circ 0=0$. Moreover, we have%
\begin{equation*}
\begin{array}{ccccc}
h+h_{1} & = & \pi \circ \iota \circ h+\pi \circ h_{1}^{\ast } & = & \pi
\circ (\iota \circ h+h_{1}^{\ast }) \\
& = & \pi \circ (\iota \circ h^{\prime }+h_{2}^{\ast }) & = & \pi \circ
\iota \circ h^{\prime }+\pi \circ h_{2}^{\ast } \\
& = & h^{\prime }+h_{2}.\blacksquare &  &
\end{array}%
\end{equation*}
\end{punto}

\begin{prop}
\label{prod-inj}Let $M$ be a left $S$-semimodule and $\{J_{\lambda
}\}_{\lambda \in \Lambda }$ be a collection of left $S$-semimodules. Then $%
\prod\limits_{\lambda \in \Lambda }J_{\lambda }$ is ($M$)-$e$-injective if
and only if $J_{\lambda }$ is $M$-$e$-injective for every $\lambda \in
\Lambda .$
\end{prop}

\begin{Beweis}
Let $J:=\prod\limits_{\lambda \in \Lambda }J_{\lambda }$ and, for each $%
\lambda \in \Lambda ,$ let $\iota _{\lambda }:J_{\lambda }\longrightarrow J$
and $\pi _{\lambda }:J\longrightarrow J_{\lambda }$ be the canonical $S$%
-linear maps.

($\Longrightarrow $) For each $\lambda \in \Lambda ,$ we have $\pi _{\lambda
}\circ \iota _{\lambda }=id_{J_{\lambda }},$ \emph{i.e.} $J_{\lambda }$ is a
retract of $J.$ The result follows from Lemma \ref{ret-inj}.

($\Longleftarrow $) Assume that $J_{\lambda }$ is $M$-$e$-injective for
every $\lambda \in \Lambda .$ Let $f:L\rightarrow M$ be a normal
monomorphism and $g:L\rightarrow J$ an $S$-linear map.%
\begin{equation*}
\xymatrix{0 \ar[r] & L \ar[r]^{f} \ar[d]_{g} & M
\ar@{-->}[ldd]^{h^*_\lambda}\\ & J \ar[d]_{\pi_\lambda} \\ & J_\lambda }
\end{equation*}%
Since $J_{\lambda }$ is $M$-$e$-injective for each $\lambda \in \Lambda ,$
there is an $S$-linear map $h_{\lambda }^{\ast }:M\rightarrow J_{\lambda }$
such that $h_{\lambda }^{\ast }\circ f=\pi _{\lambda }\circ g.$ By the \emph{%
Universal Property of Direct Products}, there exists an $S$-linear map%
\begin{equation*}
h:M\longrightarrow J,\text{ }m\mapsto \prod\limits_{\lambda \in \Lambda }({%
\iota _{\lambda }\circ h_{\lambda }^{\ast })(m).}
\end{equation*}%
Notice that for every $l\in L,$ we have
\begin{equation*}
(h\circ f)(l)=\prod\limits_{\lambda \in \Lambda }({\iota _{\lambda }\circ
h_{\lambda }^{\ast })(f(l))}=\prod\limits_{\lambda \in \Lambda }({\iota
_{\lambda }\circ \pi _{\lambda })(g(l))}=g(l).
\end{equation*}

Suppose that there exists an $S$-linear map $h^{\prime }:M\rightarrow J$
such that $h^{\prime }\circ f=g.$ It follows that $\pi _{\lambda }\circ
h^{\prime }\circ f=\pi _{\lambda }\circ g$ for every $\lambda \in \Lambda $.
Since $J_{\lambda }$ is $M$-$e$-injective, there exist $S$-linear maps $%
h_{1_{\lambda }}^{\ast },h_{2_{\lambda }}^{\ast }:M\rightarrow J$ such that $%
h_{1_{\lambda }}^{\ast }\circ f=0=h_{2_{\lambda }}^{\ast }\circ f$ and $%
h_{\lambda }^{\ast }+h_{1_{\lambda }}^{\ast }=\pi _{\lambda }\circ h^{\prime
}+h_{2_{\lambda }}^{\ast }.$%
\begin{equation*}
\xymatrix{0 \ar[r] & L \ar[d]_{g} \ar[r]^{f} & M
\ar@<-.5ex>[r]_{h^*_{2_\lambda}} \ar@<.5ex>[r]^{h^*_{1_\lambda}} \ar[ld]^{h}
\ar@/^-0.5pc/[ld]_{h'} \ar@/^0.5pc/[ldd]^{h^*_\lambda} & J_\lambda\\ & J
\ar[d]_{\pi_\lambda} \\ & J_\lambda }
\end{equation*}%
For $i=1,2,$ there exists by the \emph{Universal Property of Direct Products
}an $S$-linear map%
\begin{equation*}
h_{i}:M\longrightarrow J,\text{ }m\mapsto \prod\limits_{\lambda \in \Lambda
}(\iota _{\lambda }\circ h_{i_{\lambda }}^{\ast })(m).
\end{equation*}%
For $i=1,2$ and every $l\in L$ we have
\begin{equation*}
(h_{i}\circ f)(l)=\prod\limits_{\lambda \in \Lambda }(\iota _{\lambda }\circ
h_{i_{\lambda }}^{\ast })(f(l))=\prod\limits_{\lambda \in \Lambda }\iota
_{\lambda }(0)=0.
\end{equation*}%
Moreover, we have for every $m\in M:$%
\begin{equation*}
\begin{tabular}{lll}
$(h+h_{1})(m)$ & $=$ & $\prod\limits_{\lambda \in \Lambda }(\iota _{\lambda
}\circ \pi _{\lambda }\circ h)(m)+\prod\limits_{\lambda \in \Lambda }(\iota
_{\lambda }\circ h_{1_{\lambda }}^{\ast })(m)$ \\
& $=$ & $\prod\limits_{\lambda \in \Lambda }(\iota _{\lambda }\circ (\pi
_{\lambda }\circ h+h_{1_{\lambda }}^{\ast }))(m)$ \\
& $=$ & $\prod\limits_{\lambda \in \Lambda }(\iota _{\lambda }\circ
(h_{\lambda }^{\ast }+h_{1_{\lambda }}^{\ast }))(m)$ \\
& $=$ & $\prod\limits_{\lambda \in \Lambda }(\iota _{\lambda }\circ (\pi
_{\lambda }\circ h^{\prime }+h_{2_{\lambda }}^{\ast }))(m)$ \\
& $=$ & $\prod\limits_{\lambda \in \Lambda }(h^{\prime }+\iota _{\lambda
}\circ h_{1_{\lambda }}^{\ast })(m)$ \\
& $=$ & $(h^{\prime }+h_{2})(m).\blacksquare $%
\end{tabular}%
\end{equation*}
\end{Beweis}

\begin{lem}
\label{seq-inj}Let%
\begin{equation*}
0\longrightarrow L\overset{p}{\longrightarrow }M\overset{q}{\longrightarrow }%
N\rightarrow 0
\end{equation*}%
be a short exact sequence of left $S$-semimodules. If a left $S$-semimodule $%
J$ is $M$-$e$-injective, then $J$ is $L$-$e$-injective and $N$-$e$-injective.
\end{lem}

\begin{Beweis}
Let $J$ be a left $S$-semimodule.

\textbf{Step I:\ }$J$ is $L$-$e$-injective.

Let $f:K\rightarrow L$ be a normal monomorphism and $g:K\rightarrow J$ an $S$%
-linear map. Clearly, $p\circ f$ is a normal monomorphism.%
\begin{equation*}
\xymatrix{0 \ar[r] & K \ar[r]^{f} \ar[d]_{g} & L \ar[r]^{p} & M
\ar@{-->}[lld]^{h^*}\\ & J }
\end{equation*}%
Since $J$ is $M$-$e$-injective, there exists an $S$-linear map $h^{\ast
}:M\rightarrow J$ such that $h^{\ast }\circ p\circ f=g.$ Consider $%
h:=h^{\ast }\circ p:L\longrightarrow J.$ Then $h\circ f=h^{\ast }\circ
p\circ f=g.$

Suppose now that $h^{\prime }:L\rightarrow J$ is an $S$-linear map such that
$h^{\prime }\circ f=g$. Since $p:L\longrightarrow M$ is a normal
monomorphism and $J$ is $M$-$e$-injective, there exists an $S$-linear map $%
\tilde{h}:M\rightarrow J$ such that $\tilde{h}\circ p=h^{\prime }$. Since $%
\tilde{h}\circ p=h^{\ast }\circ p\circ f=g$, there exist $S$-linear maps $%
h_{1}^{\ast },h_{2}^{\ast }:M\rightarrow J$ such that $h_{1}^{\ast }\circ
p\circ f=0=h_{2}^{\ast }\circ p\circ f$ and $h^{\ast }+h_{1}^{\ast }=\tilde{h%
}+h_{2}^{\ast }$. Considering $h_{1}:=h_{1}^{\ast }\circ p$ and $%
h_{2}:=h_{2}^{\ast }\circ p$, we have $h_{1}\circ f=h_{1}^{\ast }\circ
p\circ f=0=h_{2}^{\ast }\circ p\circ f=h_{2}\circ f$ and
\begin{equation*}
\begin{array}{ccccc}
h+h_{1} & = & h^{\ast }\circ p+h_{1}^{\ast }\circ p & = & (h^{\ast
}+h_{1}^{\ast })\circ p \\
& = & (\tilde{h}+h_{2}^{\ast })\circ p & = & \tilde{h}\circ p+h_{2}^{\ast
}\circ p \\
& = & h^{\prime }+h_{2}. &  &
\end{array}%
\end{equation*}

\textbf{Step II:\ }$J$ is $N$-$e$-injective.

Let $f:K\rightarrow N$ be a normal monomorphism and $g:K\rightarrow J$ an $S$%
-linear map.%
\begin{equation*}
\xymatrix{& U \ar[r]^{f'} \ar[d]_{q'} & M \ar[d]_{q} \ar@/^3.0
pc/@{-->}[ldd]^{h^*} \\ 0 \ar[r] & K \ar[r]^{f} \ar[d]_{g} & N \\ & J }
\end{equation*}%
Let $(U;q^{\prime },f^{\prime })$ be a pullback of $(q,f)$ (see \cite[1.7]%
{Tak1982b}). Clearly, $f^{\prime }$ is a normal $S$-monomorphism. Since $J$
is $M$-$e$-injective, there exists an $S$-linear map $h^{\ast }:M\rightarrow
J$ such that $h^{\ast }\circ f^{\prime }=g\circ q^{\prime }$. Let $n\in N.$
Since $q$ is surjective, there exists $m_{n}\in M$ such that $n=q(m_{n}).$
Define%
\begin{equation*}
h:N\rightarrow J,\text{ }n\mapsto h^{\ast }(m_{n}).
\end{equation*}%
\textbf{Claim:}$\ h$ is well-defined.

Suppose that there exists another $m\in M$ such that $q(m)=n=q(m_{n})$.
Since $q$ is $k$-normal, there exist $m_{1},m_{2}\in Ker(q)$ such that $%
m+m_{1}=m_{n}+m_{2}.$ Since $m_{1},m_{2}\in Ker(q)$, $(m_{1},0),(m_{2},0)\in
U$ and so for $i=1,2$ we have $h^{\ast }(m_{i})=(h^{\ast }\circ f^{\prime
})(m_{i},0)=(g\circ q^{\prime })(m_{i},0)=g(0)=0,$ whence $h^{\ast
}(m)=h^{\ast }(m_{n}).$ 
Thus $h$ well defined as a map. Clearly, $h$ is $S$-linear. Moreover, for
every $k\in K$ we have $f(k)=q(m_{f(k)})$ for some $m_{f(k)}\in M,$ thus $%
(m_{f(k)},k)\in U$ and it follows that%
\begin{equation*}
\begin{array}{ccccc}
(h\circ f)(k) & = & (h\circ f\circ q^{\prime })(m_{f(k)},k) & = & (h\circ
q\circ f^{\prime })(m_{f(k)},k) \\
& = & (h^{\ast }\circ f^{\prime })(m_{f(k)},k) & = & (g\circ q^{\prime
})(m_{f(k)},k) \\
& = & g(k), &  &
\end{array}%
\end{equation*}%
i.e. $h\circ f=g$.

Suppose that there exists an $S$-linear map $h^{\prime }:N\rightarrow J$
such that $h^{\prime }\circ f=g.$ Notice that $h^{\prime }\circ q\circ
f^{\prime }=h^{\prime }\circ f\circ q^{\prime }=g\circ q^{\prime }$. Since $%
J $ is $M$-$e$-injective, there exist $h_{1}^{\ast },h_{2}^{\ast
}:M\rightarrow J$ such that $h_{1}^{\ast }\circ f^{\prime }=0=h_{2}^{\ast
}\circ f^{\prime }$ and $h^{\ast }+h_{1}^{\ast }=h^{\prime }\circ
q+h_{2}^{\ast }$.

Let $n\in N.$ Since $q$ is surjective, there exists $m_{n}\in M$ such that $%
q(m_{n})=n.$ Define
\begin{equation*}
h_{1}:N\rightarrow J,\text{ }n\mapsto h_{1}^{\ast }(m_{n})\text{ and }%
h_{2}:N\rightarrow J,\text{ }n\mapsto h_{2}^{\ast }(m_{n}).
\end{equation*}%
One can prove as above that $h_{1}$ and $h_{2}$ are well-defined. It is
clear that both $h_{1}$ and $h_{2}$ are $S$-linear. Notice that for every $%
k\in K,$ we have $(m_{f(k)},k)\in U$ whence, for $i=1,2,$ we have%
\begin{equation*}
\begin{array}{ccccc}
(h_{i}\circ f)(k) & = & (h_{i}\circ f\circ q)(m_{f(k)},k) & = & (h_{i}\circ
q\circ f^{\prime })(m_{f(k)},k) \\
& = & (h_{i}^{\ast }\circ f^{\prime })(m_{f(k)},k) & = & 0.%
\end{array}%
\end{equation*}%
Moreover, for every $n\in N$, we have%
\begin{equation*}
\begin{array}{ccccc}
(h+h_{1})(n) & = & h(n)+h_{1}(n) & = & h^{\ast }(m_{n})+h_{1}^{\ast }(m_{n})
\\
& = & (h^{\ast }+h_{1}^{\ast })(m_{n}) & = & (h^{\prime }\circ q+h_{2}^{\ast
})(m_{n}) \\
& = & (h^{\prime }\circ q)(m_{n})+h_{2}^{\ast }(m_{n}) & = & h^{\prime
}(n)+h_{2}(n) \\
& = & (h^{\prime }+h_{2})(n).\blacksquare &  &
\end{array}%
\end{equation*}
\end{Beweis}

\begin{rem}
The converse of Lemma \ref{seq-inj} is not true in general as will be shown
in Example \ref{endsnotmid}.
\end{rem}

\subsection{A Counter Example\label{sub-k-not-e-inj}}

This subsection is devoted to studying the left self-injectivity of $%
S:=M_{2}(\mathbb{R}^{+}).$ We show in particular that $\mathcal{I}%
_{S}^{i}(S)=$ $_{S}\mathbf{SM}$ and that the inclusion $\mathcal{I}%
_{S}^{e}(S)\varsubsetneqq \mathcal{I}_{S}^{i}(S)$ is strict.

\begin{lem}
\label{M2R-si}The only non-trivial proper \emph{subtractive} left ideals of $%
S$ are
\begin{eqnarray*}
E_{1} &=&Span\left( \left\{ \left[ {%
\begin{array}{cc}
1 & 0 \\
0 & 0%
\end{array}%
}\right] \right\} \right) =\left\{ \left[ {%
\begin{array}{cc}
a & 0 \\
b & 0%
\end{array}%
}\right] |\text{ }a,b\in \mathbb{R}^{+}\right\} \\
E_{2} &=&Span\left\{ \left[ {%
\begin{array}{cc}
0 & 0 \\
0 & 1%
\end{array}%
}\right] \right\} =\left\{ \left[ {%
\begin{array}{cc}
0 & c \\
0 & d%
\end{array}%
}\right] |\text{ }c,d\in \mathbb{R}^{+}\right\} \\
N_{r} &=&\left\{ \left[ {%
\begin{array}{cc}
ra & a \\
rb & b%
\end{array}%
}\right] |\text{ }a,b\in \mathbb{R}^{+}\right\} ,\text{ }r\in \mathbb{R}%
^{+}\backslash \{0\}.
\end{eqnarray*}
\end{lem}

\begin{Beweis}
We give the proof is three steps.

\textbf{Step I:}$\ E_{1},$ $E_{2}$ and $N_{r}$ are \emph{subtractive} left
ideals of $S.$

Notice that $E_{1}$ is a left ideal of $S,$ since for every $%
a,b,c,d,p,q,r,s\in \mathbb{R}^{+}$ we have
\begin{equation*}
\left[ {%
\begin{array}{cc}
p & q \\
r & s%
\end{array}%
}\right] \left[ {%
\begin{array}{cc}
a & 0 \\
b & 0%
\end{array}%
}\right] +\left[ {%
\begin{array}{cc}
c & 0 \\
d & 0%
\end{array}%
}\right] =\left[ {%
\begin{array}{cc}
pa+qb+c & 0 \\
ra+sb+d & 0%
\end{array}%
}\right] \in E_{1}.
\end{equation*}%
Moreover, $E_{1}$ is subtractive since%
\begin{equation*}
\left[ {%
\begin{array}{cc}
p & q \\
r & s%
\end{array}%
}\right] +\left[ {%
\begin{array}{cc}
a & 0 \\
b & 0%
\end{array}%
}\right] =\left[ {%
\begin{array}{cc}
c & 0 \\
d & 0%
\end{array}%
}\right]
\end{equation*}%
implies $q=0=s,$ i.e. $\left[ {%
\begin{array}{cc}
p & q \\
r & s%
\end{array}%
}\right] \in E_{1}$. Similarly, we have $E_{2}$ is a subtractive ideal.

For any non-zero $r\in \mathbb{R}^{+}$, $N_{r}$ is a left ideal since
\begin{equation*}
\left[ {%
\begin{array}{cc}
k & l \\
m & n%
\end{array}%
}\right] \left[ {%
\begin{array}{cc}
ra & a \\
rb & b%
\end{array}%
}\right] +\left[ {%
\begin{array}{cc}
rc & c \\
rd & d%
\end{array}%
}\right] =\left[ {%
\begin{array}{cc}
r(ka+lb+c) & ka+lb+c \\
r(ma+nb+d) & ma+nb+d%
\end{array}%
}\right] \in N_{r}
\end{equation*}%
for all $a,b,c,d,k,l,m,n\in \mathbb{R}^{+}$. Moreover, $N_{r}$ is
subtractive since
\begin{equation*}
\left[ {%
\begin{array}{cc}
k & l \\
m & n%
\end{array}%
}\right] +\left[ {%
\begin{array}{cc}
ra & a \\
rb & b%
\end{array}%
}\right] =\left[ {%
\begin{array}{cc}
rc & c \\
rd & d%
\end{array}%
}\right]
\end{equation*}%
implies $c=a+k/r=a+l$ and $d=b+m/r=b+n$, whence $k=rl$ and $m=rn$, i.e. $%
\left[ {%
\begin{array}{cc}
k & l \\
m & n%
\end{array}%
}\right] \in N_{r}$.

\textbf{Step II:}$\ $The only subtractive left ideal containing $E_{1},$ $%
E_{2}$ or $N_{r}$ for some $r\in \mathbb{R}^{+}\backslash \{0\}$)$\}$
strictly is $I=S.$

Let $I$ be a subtractive left ideal of $M_{2}(\mathbb{R}^{+}).$

\textbf{Case 1:\ }$E_{1}\subsetneqq I$.

In this case, there exists $\left[ {%
\begin{array}{cc}
p & q \\
r & s%
\end{array}%
}\right] \in I$ such that $q\neq 0$ or $s\neq 0$, which implies $\left[ {%
\begin{array}{cc}
0 & q \\
0 & s%
\end{array}%
}\right] \in I$ as $\left[ {%
\begin{array}{cc}
p & 0 \\
r & 0%
\end{array}%
}\right] \in I$ and%
\begin{equation*}
\left[ {%
\begin{array}{cc}
p & 0 \\
r & 0%
\end{array}%
}\right] +\left[ {%
\begin{array}{cc}
0 & q \\
0 & s%
\end{array}%
}\right] =\left[ {%
\begin{array}{cc}
p & q \\
r & s%
\end{array}%
}\right] \in I
\end{equation*}%
If $q\neq 0$, then%
\begin{equation*}
\left[ {%
\begin{array}{cc}
0 & 0 \\
0 & 1%
\end{array}%
}\right] =\left[ {%
\begin{array}{cc}
0 & 0 \\
1/q & 0%
\end{array}%
}\right] \left[ {%
\begin{array}{cc}
0 & q \\
0 & s%
\end{array}%
}\right] \in I.
\end{equation*}%
If $s\neq 0$, then
\begin{equation*}
\left[ {%
\begin{array}{cc}
0 & 0 \\
0 & 1%
\end{array}%
}\right] =\left[ {%
\begin{array}{cc}
0 & 0 \\
0 & 1/s%
\end{array}%
}\right] \left[ {%
\begin{array}{cc}
0 & q \\
0 & s%
\end{array}%
}\right] \in I.
\end{equation*}%
Either way $\left[ {%
\begin{array}{cc}
0 & 0 \\
0 & 1%
\end{array}%
}\right] \in I$, whence $E_{2}\subseteq I$ and $I=S$.

\textbf{Case 2:\ }$E_{2}\subsetneqq I$. One can show, in a was similar to
that of Case 1, that $I=S$.

\textbf{Case 3:\ } $N_{r}\subsetneqq I$ for some $r\in \mathbb{R}%
^{+}\backslash \{0\}$.

In this case, there exists some $\left[ {%
\begin{array}{cc}
k & l \\
m & n%
\end{array}%
}\right] \in I$ with $k\neq rl$ or $m\neq rn$. Assume, without loss of
generality, that $k<rl$. Then $k+p=rl$ for some $p\in \mathbb{R}%
^{+}\backslash \{0\}$. Thus $\left[ {%
\begin{array}{cc}
p & 0 \\
q & 0%
\end{array}%
}\right] \in I$ or $\left[ {%
\begin{array}{cc}
p & 0 \\
0 & q%
\end{array}%
}\right] \in I$ for some $q\in \mathbb{R}^{+}$ as%
\begin{equation*}
\left[ {%
\begin{array}{cc}
p & 0 \\
q & 0%
\end{array}%
}\right] +\left[ {%
\begin{array}{cc}
k & l \\
m & n%
\end{array}%
}\right] =\left[ {%
\begin{array}{cc}
rl & l \\
rn & n%
\end{array}%
}\right] \in I
\end{equation*}%
or
\begin{equation*}
\left[ {%
\begin{array}{cc}
p & 0 \\
0 & q%
\end{array}%
}\right] +\left[ {%
\begin{array}{cc}
k & l \\
m & n%
\end{array}%
}\right] =\left[ {%
\begin{array}{cc}
rl & l \\
m & m/r%
\end{array}%
}\right] \in I.
\end{equation*}%
Thus
\begin{equation*}
\left[ {%
\begin{array}{cc}
1 & 0 \\
0 & 0%
\end{array}%
}\right] =\left[ {%
\begin{array}{cc}
1/p & 0 \\
0 & 0%
\end{array}%
}\right] \left[ {%
\begin{array}{cc}
p & 0 \\
q & 0%
\end{array}%
}\right] \in I
\end{equation*}%
or
\begin{equation*}
\left[ {%
\begin{array}{cc}
1 & 0 \\
0 & 0%
\end{array}%
}\right] =\left[ {%
\begin{array}{cc}
1/p & 0 \\
0 & 0%
\end{array}%
}\right] \left[ {%
\begin{array}{cc}
p & 0 \\
0 & q%
\end{array}%
}\right] \in I.
\end{equation*}%
Either way we have $\left[ {%
\begin{array}{cc}
1 & 0 \\
0 & 0%
\end{array}%
}\right] \in I$, whence $E_{1}\subsetneqq I$ and so $I=S$.

\textbf{Step III.} Let $I$ be any non-zero subtractive left ideal of $S.$
Then $\left[ {%
\begin{array}{cc}
k & l \\
m & n%
\end{array}%
}\right] \in I\backslash \{0\}$ for some $k,l,m,n\in \mathbb{R}^{+}$.

\textbf{Case 1:} $k\neq 0.$ In this case, we have%
\begin{equation*}
\left[ {%
\begin{array}{cc}
1/k & 0 \\
0 & 0%
\end{array}%
}\right] \left[ {%
\begin{array}{cc}
k & l \\
m & n%
\end{array}%
}\right] =\left[ {%
\begin{array}{cc}
1 & l/k \\
0 & 0%
\end{array}%
}\right] \in I.
\end{equation*}%
If $l=0,$ then $\left[ {%
\begin{array}{cc}
1 & 0 \\
0 & 0%
\end{array}%
}\right] \in I,$ whence $E_{1}\subseteq I.$ Otherwise, $\left[ {%
\begin{array}{cc}
k/l & 1 \\
0 & 0%
\end{array}%
}\right] \in I,$ whence $N_{k/l}\subseteq I.$ In either case, it follows by
Step II that $I\in \{E_{1},N_{k/l},S\}.$

\textbf{Case 2:} $l\neq 0$. In this case, we have%
\begin{equation*}
\left[ {%
\begin{array}{cc}
0 & 0 \\
1/l & 0%
\end{array}%
}\right] \left[ {%
\begin{array}{cc}
k & l \\
m & n%
\end{array}%
}\right] =\left[ {%
\begin{array}{cc}
0 & 0 \\
k/l & 1%
\end{array}%
}\right] \in I.
\end{equation*}%
If $k=0,$ then $\left[ {%
\begin{array}{cc}
0 & 0 \\
0 & 1%
\end{array}%
}\right] \in I,$ whence $E_{2}\subseteq I.$ Otherwise $N_{k/l}\subseteq I.$
In either case, it follows by Step II that $I\in \{E_{2},N_{k/l},S\}.$

\textbf{Case 3:} $m\neq 0$. In this case, we have%
\begin{equation*}
\left[ {%
\begin{array}{cc}
0 & 1/m \\
0 & 0%
\end{array}%
}\right] \left[ {%
\begin{array}{cc}
k & l \\
m & n%
\end{array}%
}\right] =\left[ {%
\begin{array}{cc}
1 & n/m \\
0 & 0%
\end{array}%
}\right] \in I.
\end{equation*}%
If $n=0,$ then $\left[ {%
\begin{array}{cc}
1 & 0 \\
0 & 0%
\end{array}%
}\right] \in I,$ whence $E_{1}\subseteq I.$ Otherwise, $\left[ {%
\begin{array}{cc}
m/n & 1 \\
0 & 0%
\end{array}%
}\right] \in I,$ whence $N_{m/n}\subseteq I.$ In either case, it follows by
Step II that $I\in \{E_{1},N_{m/n},S\}.$

\textbf{Case 4}:\ $n\neq 0$. In this case, we have%
\begin{equation*}
\left[ {%
\begin{array}{cc}
0 & 0 \\
0 & 1/n%
\end{array}%
}\right] \left[ {%
\begin{array}{cc}
k & l \\
m & n%
\end{array}%
}\right] =\left[ {%
\begin{array}{cc}
0 & 0 \\
m/n & 1%
\end{array}%
}\right] .
\end{equation*}%
If $m=0,$ then $\left[ {%
\begin{array}{cc}
0 & 0 \\
0 & 1%
\end{array}%
}\right] \in I,$ whence $E_{2}\subseteq I.$ Otherwise, $N_{m/n}\subseteq I.$
In either case, it follows by Step II that $I\in
\{E_{2},N_{m/n},S\}.\blacksquare $
\end{Beweis}

\begin{lem}
\label{M2R-all-i-inj}Every left $S$-semimodule is $S$-$i$-injective.
\end{lem}

\begin{Beweis}
Let $M$ be a left $S$-semimodule, $f:N\rightarrow S$ a normal $S$%
-monomorphism, and $g:N\rightarrow M$ an $S$-linear map. Then $f(N)$ is a
subtractive left ideal of $S$, whence $f(N)\in \{0,E_{1},E_{2},S\}$ or $%
f(N)=N_{r}$ for some $r\in \mathbb{R}^{+}\backslash \{0\}$.

\textbf{Case I:}\ $f(N)=0.$ In this case, choose $h=0:S\rightarrow M$.
Clearly, $g=h\circ f$.

\textbf{Case II:}\ $f(N)=S$. In this case, $f$ is an $S$-isomorphism. Choose
$h=g\circ f^{-1},$ whence $g=h\circ f$.

\textbf{Case III:}\ $f(N)=E_{1}$. In this case, there exists a unique $%
n_{0}\in N$ such that
\begin{equation*}
f(n_{0})=\left[ {%
\begin{array}{cc}
1 & 0 \\
0 & 0%
\end{array}%
}\right] .
\end{equation*}%
Consider the $S$-linear map%
\begin{equation*}
h:S\rightarrow M,\text{ }\left[ {%
\begin{array}{cc}
p & q \\
r & s%
\end{array}%
}\right] \longmapsto \left[ {%
\begin{array}{cc}
p & q \\
r & s%
\end{array}%
}\right] g(n_{0}).
\end{equation*}%
Let $n\in N$. It follows that%
\begin{equation*}
f(n)=\left[ {%
\begin{array}{cc}
a & 0 \\
b & 0%
\end{array}%
}\right] =\left[ {%
\begin{array}{cc}
a & 0 \\
b & 0%
\end{array}%
}\right] f(n_{0})=f\left( \left[ {%
\begin{array}{cc}
a & 0 \\
b & 0%
\end{array}%
}\right] n_{0}\right)
\end{equation*}%
for some $a,b\in \mathbb{R}^{+}$. Since $f$ is injective, $n=\left[ {%
\begin{array}{cc}
a & 0 \\
b & 0%
\end{array}%
}\right] n_{0}$. It follows that%
\begin{equation*}
\begin{array}{ccccc}
(h\circ f)(n) & = & h(f(n)) & = & h(\left[ {%
\begin{array}{cc}
a & 0 \\
b & 0%
\end{array}%
}\right] ) \\
& = & \left[ {%
\begin{array}{cc}
a & 0 \\
b & 0%
\end{array}%
}\right] g(n_{0}) & = & g\left( \left[ {%
\begin{array}{cc}
a & 0 \\
b & 0%
\end{array}%
}\right] n_{0}\right) \\
& = & g(n). &  &
\end{array}%
\end{equation*}

\textbf{Case IV:}\ $f(N)=E_{2}.$ The proof is similar to Case III.

\textbf{Case V:} $f(N)=N_{r}$ for some $r\in \mathbb{R}^{+}\backslash \{0\}.$
In this case, there exists a unique $n_{0}\in N$ such that%
\begin{equation*}
f(n_{0})=\left[ {%
\begin{array}{cc}
1 & 1/r \\
0 & 0%
\end{array}%
}\right] .
\end{equation*}%
Define
\begin{equation*}
h:S\rightarrow M,\text{ }\left[ {%
\begin{array}{cc}
j & k \\
l & m%
\end{array}%
}\right] \longmapsto \left[ {%
\begin{array}{cc}
j & k \\
l & m%
\end{array}%
}\right] g(n_{0}).
\end{equation*}%
For every $n\in N$, we have%
\begin{equation*}
f(n)=\left[ {%
\begin{array}{cc}
ra & a \\
rb & b%
\end{array}%
}\right] =\left[ {%
\begin{array}{cc}
ra & a \\
rb & b%
\end{array}%
}\right] f(n_{0})=f\left( \left[ {%
\begin{array}{cc}
ra & a \\
rb & b%
\end{array}%
}\right] n_{0}\right)
\end{equation*}%
for some $a,b\in \mathbb{R}^{+}$. Since $f$ is injective, $n=\left[ {%
\begin{array}{cc}
ra & a \\
rb & b%
\end{array}%
}\right] n_{0}$ and so%
\begin{equation*}
\begin{array}{ccccc}
(h\circ f)(n) & = & h(f(n)) & = & h\left( \left[ {%
\begin{array}{cc}
ra & a \\
rb & b%
\end{array}%
}\right] \right) \\
& = & \left[ {%
\begin{array}{cc}
ra & a \\
rb & b%
\end{array}%
}\right] g(n_{0}) & = & g\left( \left[ {%
\begin{array}{cc}
ra & a \\
rb & b%
\end{array}%
}\right] n_{0}\right) \\
& = & g(n).\blacksquare &  &
\end{array}%
\end{equation*}
\end{Beweis}

We are now ready to provide an example of an $S$-$i$-injective semimodule
which is not $S$-$e$-injective.

\begin{ex}
\label{k-inj-not-e-ink}The left $S$-semimodule%
\begin{equation}
N_{1}=\left\{ \left[ {%
\begin{array}{cc}
a & a \\
b & b%
\end{array}%
}\right] |\text{ }a,b\in \mathbb{R}^{+}\right\}  \label{N1}
\end{equation}%
is $S$-$i$-injective but not $S$-$e$-injective.
\end{ex}

\begin{Beweis}
Let $\iota :N_{1}\rightarrow S$ be an embedding and $id:N_{1}\rightarrow
N_{1}$ be the identity map. Since $N_{1}$ is subtractive, $\iota $ is a
normal $S$-monomorphism. Let $h_{1},h_{2}:S\rightarrow N_{1}$ with
\begin{equation*}
h_{1}\left( \left[ {%
\begin{array}{cc}
p & q \\
r & s%
\end{array}%
}\right] \right) =\left[ {%
\begin{array}{cc}
p & p \\
r & r%
\end{array}%
}\right] \text{ and }h_{2}\left( \left[ {%
\begin{array}{cc}
p & q \\
r & s%
\end{array}%
}\right] \right) =\left[ {%
\begin{array}{cc}
q & q \\
s & s%
\end{array}%
}\right] .
\end{equation*}%
Then%
\begin{equation*}
\begin{array}{ccc}
(h_{1}\circ \iota )\left( \left[ {%
\begin{array}{cc}
a & a \\
b & b%
\end{array}%
}\right] \right) & = & h_{1}\left( \left[ {%
\begin{array}{cc}
a & a \\
b & b%
\end{array}%
}\right] \right) \\
& = & \left[ {%
\begin{array}{cc}
a & a \\
b & b%
\end{array}%
}\right] \\
& = & id\left( \left[ {%
\begin{array}{cc}
a & a \\
b & b%
\end{array}%
}\right] \right)%
\end{array}%
\end{equation*}
and%
\begin{equation*}
\begin{array}{ccc}
(h_{2}\circ \iota )\left( \left[ {%
\begin{array}{cc}
a & a \\
b & b%
\end{array}%
}\right] \right) & = & h_{2}\left( \left[ {%
\begin{array}{cc}
a & a \\
b & b%
\end{array}%
}\right] \right) \\
& = & \left[ {%
\begin{array}{cc}
a & a \\
b & b%
\end{array}%
}\right] \\
& = & id\left( \left[ {%
\begin{array}{cc}
a & a \\
b & b%
\end{array}%
}\right] \right) .%
\end{array}%
\end{equation*}

Suppose that there exist $k_{1},k_{2}:S\rightarrow N_{1}$ such that $%
k_{1}\circ \iota =0=k_{2}\circ \iota $ and $h_{1}+k_{1}=h_{2}+k_{2}$. Write
\begin{equation*}
k_{1}\left( \left[ {%
\begin{array}{cc}
1 & 0 \\
0 & 1%
\end{array}%
}\right] \right) =\left[ {%
\begin{array}{cc}
l & m \\
n & o%
\end{array}%
}\right] \text{ and }k_{2}\left( \left[ {%
\begin{array}{cc}
1 & 0 \\
0 & 1%
\end{array}%
}\right] \right) =\left[ {%
\begin{array}{cc}
p & q \\
r & s%
\end{array}%
}\right]
\end{equation*}%
for some $k,l,m,n,o,p,q,r,s\in \mathbb{R}^{+}$. Then
\begin{eqnarray*}
k_{1}\left( \left[ {%
\begin{array}{cc}
a & b \\
c & d%
\end{array}%
}\right] \right) &=&\left[ {%
\begin{array}{cc}
a & b \\
c & d%
\end{array}%
}\right] \left[ {%
\begin{array}{cc}
l & m \\
n & o%
\end{array}%
}\right] \\
k_{2}\left( \left[ {%
\begin{array}{cc}
a & b \\
c & d%
\end{array}%
}\right] \right) &=&\left[ {%
\begin{array}{cc}
a & b \\
c & d%
\end{array}%
}\right] \left[ {%
\begin{array}{cc}
p & q \\
r & s%
\end{array}%
}\right]
\end{eqnarray*}%
for every $a,b,c,d\in \mathbb{R}^{+}$. It follows that%
\begin{equation*}
0=(k_{1}\circ \iota )\left( \left[ {%
\begin{array}{cc}
1 & 1 \\
0 & 0%
\end{array}%
}\right] \right) =\left[ {%
\begin{array}{cc}
1 & 1 \\
0 & 0%
\end{array}%
}\right] \left[ {%
\begin{array}{cc}
l & m \\
n & o%
\end{array}%
}\right] =\left[ {%
\begin{array}{cc}
l+n & m+o \\
0 & 0%
\end{array}%
}\right] ,
\end{equation*}%
which implies that $l=m=n=o=0$ as $0$ is the only element of $\mathbb{R}^{+}$
which has additive inverse. So,%
\begin{equation*}
\begin{array}{ccccc}
0 & = & \left( k_{2}\circ \iota \right) \left( \left[ {%
\begin{array}{cc}
1 & 1 \\
0 & 0%
\end{array}%
}\right] \right) & = & k_{2}\left( \left[ {%
\begin{array}{cc}
1 & 1 \\
0 & 0%
\end{array}%
}\right] \right) \\
& = & \left[ {%
\begin{array}{cc}
1 & 1 \\
0 & 0%
\end{array}%
}\right] \left[ {%
\begin{array}{cc}
p & q \\
r & s%
\end{array}%
}\right] & = & \left[ {%
\begin{array}{cc}
p+r & q+s \\
0 & 0%
\end{array}%
}\right] ,%
\end{array}%
\end{equation*}%
which implies that $p=q=r=s=0$ as $0$ is the only element of $\mathbb{R}^{+}$
which has additive inverse. Thus $k_{1}=0=k_{2}$, a contradiction with $%
h_{1}+k_{1}=h_{2}+k_{2}$ as $h_{1}\neq h_{2}$. Hence, $N_{1}$ is not $S$-$e$%
-injective.$\blacksquare $
\end{Beweis}

The following example shows that the converse of Lemma \ref{seq-inj} is not
true in general.

\begin{ex}
\label{endsnotmid}Consider the short exact sequence%
\begin{equation*}
0\rightarrow E_{1}\overset{\iota _{E_{1}}}{\longrightarrow }S\overset{\pi
_{E_{2}}}{\longrightarrow }E_{2}\rightarrow 0
\end{equation*}%
of left $S$-semimodules. Then $N_{1}$ is $E_{1}$-$e$-injective and $E_{2}$-$%
e $-injective but not $S$-$e$-injective.
\end{ex}

\begin{Beweis}
Let $f:M\rightarrow E_{1},g:M\rightarrow N_{1}$ be $S$-linear maps where $f$
is a normal monomorphism. If $f=0,$ then we are done. If $f\neq 0.$ then $f$
is an isomorphism as $E_{1}$ is ideal-simple. Define $h=g\circ f^{-1}$. Then
$h\circ f=(g\circ f^{-1})\circ f=g$. If $h^{\prime }:E_{1}\rightarrow N_{1}$
is an $S$-linear map satisfies $h^{\prime }\circ f=g$, then $h^{\prime
}=h^{\prime }\circ (f\circ f^{-1})=g\circ f^{-1}=h$. Hence $N_{1}$ is $E_{1}$%
-$e$-injective. Similarly, $N_{1}$ is $E_{2}$-$e$-injective. However, $N_{1}$
is not $S$-$e$-injective as shown in Example \ref{k-inj-not-e-ink}.
\end{Beweis}

\section{The Embedding Problem}

It is well-known that the category of left (right) modules over a ring $R$
has enough injectives, \emph{i.e. }every left (right) $R$-module $M$ can be
embedded in an injective left (right) $R$-module (e.g., $E(M),$ the \emph{%
injective hull} of $M$). This is true only for the left (right) semimodules
over some special semirings which are not rings (e.g., the \emph{additively
idempotent semirings} \cite{Wan1994}, \cite[Corollary 17.34]{Gol1999}). In
\cite{Ili2016}, Il'in conjectured that a semiring $S$ has the property that
every left (right)\ $S$-semimodule has an \emph{injective hull} if and only
if $S$ is \emph{additively regular} (i.e. for every $a\in S,$ there exists
some $x\in S$ such that $a+x+a=a$). In fact, the situation over some
semirings can be extremely bad:

\begin{lem}
\label{N0}If $S$ is an entire, cancellative, zerosumfree semiring, then the
only injective left $S$-semimodule is $\{0\}$ \emph{(}cf., \emph{\cite[%
Proposition 17.21]{Gol1999})}.
\end{lem}

\begin{ex}
The category of commutative monoids (i.e., $\mathbb{Z}^{+}$-semimodules) has
no non-zero injective objects.
\end{ex}

Another significant difference is that Baer's Criterion (a left module $M$
over a ring $R$ is injective if $M$ is $R$-injective) is not valid for
semimodules over arbitrary semirings (which are not rings).

\begin{defn}
Let $S$ be a semiring. We say that the category $_{S}\mathbf{SM}$ has
\textbf{enough injectives} (resp. \textbf{enough }$e$\textbf{-injectives},
\textbf{enough }$i$\textbf{-injectives}), if every left $S$-semimodule can
be embedded in an injective (resp. $e$-injective, $i$-injective) left $S$%
-semimodule.
\end{defn}

\begin{lem}
\emph{(\cite[Theorem 3]{Ili2008})} If $_{S}S$ satisfies the Baer's criterion
and $_{S}\mathbf{SM}$ has enough injectives, then $S$ is a ring.

\begin{prop}
\label{enj->e}Let $S$ be a semiring. If $_{S}\mathbf{SM}$ has enough $e$%
-injectives, then every injective left $S$-semimodule is $e$-injective.

\begin{Beweis}
Let $I$ be an injective left $S$-semimodule. By assumption, there is an
embedding $I\overset{\iota }{\hookrightarrow }E,$ where $_{S}E$ is $e$%
-injective. Since $_{S}I$ is injective, there exists and $S$-linear map $\pi
:E\longrightarrow I$ such that $\pi \circ \iota =id_{I}.$ It follows that $%
_{S}I$ is a \emph{retract} of an $e$-injective left $S$-semimodule, whence $%
e $-injective by Proposition \ref{ret-inj}.$\blacksquare $
\end{Beweis}
\end{prop}
\end{lem}

\begin{prop}
\label{idemp-embed-e-inj}\emph{(}compare with \emph{\cite[Theorem 4.5]%
{AIKN2018})}\ Let $S$ be an \emph{additively idempotent} semiring. Then $_{S}%
\mathbf{SM}$ has enough $e$-injectives, and every injective left $S$%
-semimodule is $e$-injective.
\end{prop}

\begin{punto}
We define a left $S$-semimodule $N$ to be \textbf{divisible}, if for every $%
s\in S,$ which is not a zero divisor, there exists for every $n\in N$ some $%
m_{n}\in N$ such that $sm_{n}=n.$ As in the case of modules over a ring,
every injective semimodule over a semiring is divisible.
\end{punto}

The proof of the following observation is similar to that in the case of
modules over rings \cite[16.6]{Wis1991}.

\begin{lem}
\label{inj->div}Every $S$-injective left $S$-semimodule is divisible.
\end{lem}

\begin{Beweis}
Let $N$ be an injective left $S$-semimodule and $n\in N.$ Let $s\in S$ be a
non zero-divisor. Claim:\ there exists $m_{n}\in N$ such that $sm_{n}=n.$
Consider the canonical embedding $0\longrightarrow Ss\overset{\iota }{%
\longrightarrow }S$ and the $S$-linear map
\begin{equation*}
h:Ss\longrightarrow N,\text{ }ts\mapsto tn.
\end{equation*}%
By our assumption, $N$ is $S$-injective, whence there exists an $S$-linear
map $g:S\longrightarrow N$ such that $g\circ \iota =h.$ Let $%
m_{n}:=g(1_{S}). $ Then we have
\begin{equation*}
n=h(s)=(g\circ \iota )(s)=g(s)=g(s\cdot 1_{S})=sg(1_{S})=sm_{n}.\blacksquare
\end{equation*}
\end{Beweis}

The converse of Lemma \ref{inj->div} is not true in general as the following
example shows.

\begin{ex}
$\mathbb{Q}$ is a divisible commutative monoid which is \emph{not} injective.
\end{ex}

\begin{punto}
Let $R$ be a ring. Every left $R$-module can be embedded in an injective
module $Hom_{\mathbb{Z}}(R,D)$, (cf., \cite[page 407, 421]{Gri2007}). For a
semiring $S$, we prove that every left $S$-semimodule can be embedded into $%
Hom_{\mathbb{Z}^{+}}(S,D)$ for some divisible commutative monoid $D$.
However, it is unknown whether $Hom_{\mathbb{Z}^{+}}(S,D)$ is necessarily $e$%
-injective.
\end{punto}

\begin{lem}
\label{x410}Every commutative monoid can be embedded \emph{subtractively} in
a divisible commutative monoid.
\end{lem}

\begin{Beweis}
Let $B$ be a commutative monoid. Then there exists a surjective morphism of
monoids $f:\mathbb{Z}^{+(\Lambda )}\rightarrow B$ for some index set $%
\Lambda .$ Let $g$ be the embedding of $\mathbb{Z}^{+(\Lambda )}$ into $%
\mathbb{Q}^{+(\Lambda )}.$ Let $(g^{\prime },f^{\prime };P)$ be a pushout of
$(f,g)$ (see \cite[Theorem 2.3, Corollary 2.4]{AN}).
\begin{equation*}
\xymatrix{\mathbb{Q}^{+(\Lambda)} \ar[r]^{f'} & P \\
{\mathbb{Z}^{+}}^{(\Lambda)} \ar[u]^{g} \ar[r]_{f} & B \ar[u]_{g'} }
\end{equation*}

Notice that $g^{\prime }$ is subtractive since $g$ is subtractive. Moreover,
the commutative monoid $P$ is divisible since for every $n\in \mathbb{Z}^{+}$
and $p\in P$ we have $(q_{\lambda })_\Lambda,(q_{\lambda
}^\prime)_\Lambda\in \mathbb{Q}^{+}$ such that $p=f^\prime((q_{\lambda
})_\Lambda)$ and $nq_{\lambda }^{\prime }=q_{\lambda }$. Thus $%
nf^\prime((q_{\lambda }^{\prime })_\Lambda)=f^\prime((nq_{\lambda }^{\prime
})_\Lambda)=f^\prime((q_{\lambda })_\Lambda).$

Let $C:=\{q\in \mathbb{Q}^{+}|0\leq q<1\}.$ Then $B\oplus C^{(\Lambda )}$ is
a commutative monoid with
\begin{equation*}
(b,(c_{\lambda }))+(b^{\prime },(c_{\lambda }^{\prime }))=(b+b^{\prime
}+f((\lfloor c_{\lambda }+c_{\lambda }^{\prime }\rfloor )_{\Lambda
}),(c_{\lambda }+c_{\lambda }^{\prime }-\lfloor c_{\lambda }+c_{\lambda
}^{\prime }\rfloor )).
\end{equation*}
\begin{equation*}
\xymatrix{& & B\oplus C^{(\Lambda)} \\ \mathbb{Q}^{+(\Lambda)} \ar[r]^{f'}
\ar@/^1 pc/[rru]^{f^*} & P \ar@{-->}[ru]^{\varphi} \\
{\mathbb{Z}^{+}}^{(\Lambda)} \ar[u]^{g} \ar[r]_{f} & B \ar[u]_{g'} \ar@/_1
pc/[ruu]_{g^*} }
\end{equation*}%
The map
\begin{equation*}
g^{\ast }:B\longrightarrow B\oplus C^{(\Lambda )},\text{ }b\mapsto (b,0)
\end{equation*}%
is a $\mathbb{Z}^{+}$-monomorphism. The map
\begin{equation*}
f^{\ast }:\mathbb{Q}^{+(\Lambda )}\longrightarrow B\oplus C^{(\Lambda )},%
\text{ }(q_{\lambda })\mapsto (f((\lfloor q_{\lambda }\rfloor )_{\Lambda
}),(q_{\lambda }-\lfloor q_{\lambda }\rfloor )_{\Lambda })
\end{equation*}%
is a $\mathbb{Z}^{+}$-homomorphism. Since $f^{\ast }\circ g=g^{\ast }\circ
f, $ there exists, by the \emph{Universal Property of Pushouts}, a unique
map $\varphi :P\rightarrow B\oplus C^{(\Lambda )}$ such that $\varphi \circ
f^{\prime }=f^{\ast }$ and $\varphi \circ g^{\prime }=g^{\ast }.$ Since $%
g^{\ast }$ is injective, $g^{\prime }$ is injective. Hence $g^{\prime
}:B\longrightarrow P$ is a normal $\mathbb{Z}^{+}$-monomorphism from $B$
into the divisible commutative monoid $P.\blacksquare $
\end{Beweis}

\begin{lem}
\label{embinj}Every left $S$-semimodule can be embedded into $Hom_{\mathbb{Z}%
^{+}}(S,D)$ for some divisible commutative monoid $D$.
\end{lem}

\begin{Beweis}
Let $M$ be a left $S$-semimodule. By Lemma \ref{x410} there exists a normal
monomorphism of commutative monoids $\mu :M\rightarrow D$ for some divisible
commutative monoid $D.$ Consider the canonical $S$-linear map%
\begin{equation*}
\epsilon :M\longrightarrow Hom_{\mathbb{Z}^{+}}(S,D),\text{ }m\mapsto
\lbrack s\mapsto \mu (sm)].
\end{equation*}%
Suppose that $\epsilon (m)=\epsilon (m^{\prime })$ for some $m,m^{\prime
}\in M.$ Then, in particular, $\epsilon (m)(1_{S})=\epsilon (m^{\prime
})(1_{S}),$ i.e. $\mu (m)=\mu (m^{\prime }).$ Since $\mu $ is injective, we
conclude that $m=m^{\prime }.\blacksquare $
\end{Beweis}

The embedding into an injective $R$-module (where $R$ is a ring) implies a
nice result in the category of $R$-modules: an $R$-module $P$ is projective
if and only if $P$ is $J$-projective for every injective $R$-module $J$ \cite%
[page 411]{Gri2007}. For semimodules, we have so far the following
implication.

\begin{prop}
\label{dtohomd}Let $\gamma :T\longrightarrow S$ be a morphism of semirings
and $M$ a left $S$-semimodule. If $_{T}A$ is $_{T}M$-$i$-injective, then $%
Hom_{T}(S,A)$ is $_{S}M$-$i$-injective.
\end{prop}

\begin{Beweis}
Let $\iota :K\rightarrow M$ be a normal $S$-monomorphism and $f:K\rightarrow
Hom_{T}(S_{S},A)$ an $S$-linear map.
\begin{equation*}
\xymatrix{0 \ar[r] & K \ar[r]^{\iota} \ar[d]^{f} & M \\ & Hom_T(S_S,A) }
\end{equation*}%
Recall the canonical \emph{isomorphism} of commutative monoids
\begin{equation*}
Hom_{S}(K,Hom_{T}(S_{S},A))\overset{\theta _{K,A}}{\simeq }Hom_{T}(K,A),%
\text{ }f\mapsto \lbrack k\mapsto f(k)(1_{S})].
\end{equation*}%
Consider the $T$-linear map $\theta _{K,A}(f):K\longrightarrow A.$%
\begin{equation*}
\xymatrix{0 \ar[r] & K \ar[r]^{\iota} \ar[d]_{\theta_{K,A}(f)} & M
\ar@{-->}[ld]^{h} \\ & A }
\end{equation*}%
Since $\iota :K\rightarrow M$ is also a normal $T$-monomorphism and $_{T}A$
is $M$-$i$-injective, there exists a $T$-linear map $h:M\longrightarrow A$
such that $h\circ \iota =\theta _{K,A}(f).$ Notice that $\theta
_{M,A}^{-1}(h):M\rightarrow Hom_{T}(S_{S},A)$ is $S$-linear and that for all
$k\in K$ and every $s\in S$ we have%
\begin{equation*}
\begin{array}{ccccc}
((\theta _{M,A}^{-1}(h)\circ \iota )(k))(s) & = & \theta
_{M,A}^{-1}(h)(s\iota (k)) & = & h(s\iota (k)) \\
& = & (h\circ \iota )(sk) & = & \theta _{K,A}(f)(sk) \\
& = & f(sk)(1_{S}) & = & (sf(k))(1_{S}) \\
& = & f(k)(1_{S}\cdot s) & = & f(k)(s).%
\end{array}%
\end{equation*}%
Hence, $Hom_{T}(S_{S},A)$ is $M$-$i$-injective as a left $S$-semimodule.$%
\blacksquare $
\end{Beweis}

The following result is a combination of Proposition \ref{dtohomd} and \cite[%
Corollary 3.5]{AIKN2018}.

\begin{cor}
\label{pres-inj}Let $\gamma :T\longrightarrow S$ be a morphism of semirings.
The functor
\begin{equation*}
Hom_{T}(S_{S},-):\text{ }_{T}\mathbf{SM}\longrightarrow \text{ }_{S}\mathbf{%
SM}
\end{equation*}%
preserves injective, $e$-injective and $i$-injective objects.
\end{cor}

\begin{lem}
\label{diviinj}Every divisible commutative monoid is $\mathbb{Z}^{+}$-$i$%
-injective.
\end{lem}

\begin{Beweis}
Let $D$ be a divisible commutative monoid, $f:I\rightarrow \mathbb{Z}^{+}$ a
normal monomorphism of commutative monoids and $g:I\rightarrow D$ a morphism
of commutative monoids. Since $f(I)$ is subtractive, $f(I)=k\mathbb{Z}^{+}$
for some $k\in \mathbb{Z}^{+}.$ Let $i_{0}\in I$ be such that $f(i_{0})=k$
and notice that $i_{0}$ is unique as $f$ is injective. By our choice, $D$ is
divisible and so there exists $d\in D$ such that $kd=g(i_{0})$. The map%
\begin{equation*}
h:\mathbb{Z}^{+}\rightarrow D,\text{ }n\mapsto nd
\end{equation*}%
is a well-defined morphism of monoids. Moreover, for every $i\in I,$ we have
$f(i)=nk$ for some $n\in \mathbb{Z}^{+}$ whence $i=ni_{0}$ as $f$ is
injective. It follows that $f(i)=f(ni_{0})=nf(i_{0})=nk,$ and so
\begin{equation*}
(h\circ f)(i)=h(nk)=h(n)k=ndk=ng(i_{0})=g(ni_{0})=g(i).
\end{equation*}%
It follows that $hf=g.\blacksquare $
\end{Beweis}

\begin{defn}
We say that a left $S$-semimodule $I$ is $c$\textbf{-injective} (resp. $c$-$%
e $-\textbf{injective}, $c$-$i$-\textbf{injective),} if $I$ is $M$-injective
(resp., $M$-$e$-injective, $M$-$i$-injective) for every cancellative left $S$%
-semimodule $M.$
\end{defn}

\begin{prop}
\label{cancl}Every divisible commutative monoid is $c$-$i$-injective.
\end{prop}

\begin{Beweis}
Let $D$ be a divisible commutative monoid, $N$ a \emph{cancellative} left $S$%
-semimodule, $f:M\rightarrow N$ a normal $\mathbb{Z}^{+}$-monomorphism and $%
g:M\rightarrow J$ a morphism of commutative monoids.%
\begin{equation*}
\xymatrix{0 \ar[r] & M \ar[r]^{f} \ar[d]^{g} & N \\ & J }
\end{equation*}%
Define
\begin{equation*}
\mathcal{S}=\{(A,\alpha ):A\leq _{\mathbb{Z}^{+}}N,\text{ }M\subseteq A,%
\text{ }\alpha :A\rightarrow J\text{ with }\alpha (m)=g(m)\text{ }\forall
\text{ }m\in M\}.
\end{equation*}%
Notice that $\mathcal{S}$ is not empty, since $(M,g)\in \mathcal{S}$. Define
an order on $\mathcal{S}$ as follows:%
\begin{equation*}
(A,\alpha )\leq (B,\beta )\Leftrightarrow A\subseteq B\text{ and }\beta
(a)=\alpha (a)\text{ }\forall a\in A.
\end{equation*}%
Let $((A_{\lambda },\alpha _{\lambda }))_{\Lambda }$ be a chain in $\mathcal{%
S}$. Set $\mathbf{A}:=\bigcup\limits_{\lambda \in \Lambda }A_{\lambda }$ and
define $\mathbf{\alpha }:A\rightarrow J$ such that, if $x\in A_{\lambda },$
then $\alpha (x)=\alpha _{\lambda }(x)$. Notice that $\alpha $ is
well-defined, thus the chain has an upper bound $(\mathbf{A},\mathbf{\alpha }%
)$. By Zorn's Lemma, $\mathcal{S}$ has a maximal element $(C,\gamma )$.

\textbf{Claim: }If $A\neq N$, then $(A,\alpha )\in \mathcal{S}$ is not
maximal.

Let $(A,\alpha )\in \mathcal{S}$ with $A\subsetneqq N.$ Choose $b\in
N\backslash A$ and set $B:=A+\mathbb{Z}^{+}b$. Notice that $L:=\{r\in
\mathbb{Z}^{+}|\text{ }rb\in A\}$ is an ideal of $\mathbb{Z}^{+}$ and
\begin{equation*}
\kappa :L\longrightarrow J,\text{ }r\mapsto \alpha (rb)
\end{equation*}%
is a morphism of monoids. By Lemma \ref{diviinj} there exists a morphism of
monoids $\chi :\mathbb{Z}^{+}\rightarrow J$ such that $\chi (r)=\alpha (rb)$
$\forall $ $r\in L$. Define%
\begin{equation*}
\beta :B\rightarrow J,\text{ }a+rb\mapsto \alpha (a)+\chi (r).
\end{equation*}%
We claim that $\beta $ is well-defined. Suppose that $a+rb=a^{\prime
}+r^{\prime }b$ for some $r\in L$ and $a\in A.$ Assume, without loss of
generality, that $r^{\prime }>r,$ whence $r^{\prime }=r+\tilde{r}$ for some $%
\tilde{r}\in \mathbb{Z}^{+}.$ It follows that $a+rb=a^{\prime }+r^{\prime
}b=a^{\prime }+rb+\tilde{r}b$, whence $a=a^{\prime }+\tilde{r}b$ as $N$ is
cancellative. It follows that%
\begin{equation*}
\begin{array}{ccccc}
\beta (a^{\prime }+r^{\prime }b) & = & \beta ((a^{\prime }+\tilde{r}b)+rb) &
= & \alpha (\tilde{r}b+a^{\prime })+\chi (r) \\
& = & \alpha (a)+\chi (r) & = & \beta (a+rb).%
\end{array}%
\end{equation*}%
Thus $\beta $ is well-defined as morphism of monoids with $\beta (a)=\alpha
(a)$ $\forall $ $a\in A$. Thus $(A,\alpha )$ is not maximal in . It follows
that there exists a morphism of monoids $h:N\longrightarrow J$ such that $%
(N,h)$ is maximal in $\mathcal{S}.$ Clearly, $h:N\rightarrow J$ such that $%
h\circ f=g.\blacksquare $
\end{Beweis}

The following result is, in some sense, a generalization of the fact
(mentioned without proof in \cite[17.35]{Gol1999}) that any \emph{%
cancellative semimodule} over semiring can be embedded in a $c$-injective
module. While $c$-$i$-injectivity is formally weaker than $c$-injectivity,
our result works for arbitrary, not necessarily cancellative, semimodules
over semirings.

\begin{thm}
\label{c-i-injective}Every left $S$-semimodule can be embedded as a \emph{%
subtractive} subsemimodule of a $c$-$i$-injective left $S$-semimodule.
\end{thm}

\begin{Beweis}
Let $M$ be a left $S$-semimodule. By Lemma \ref{embinj}, $M$ can be embedded
as a \emph{subtractive }subsemimodule of the left $S$-semimodule $Hom_{%
\mathbb{Z}^{+}}(S,D)$ for some divisible commutative monoid $D.$ Let $N$ be
a cancellative left $S$-semimodule; then $N$ is, in particular, a
cancellative commutative monoid. By Proposition \ref{cancl}, $D$ is an $N$-$%
i $-injective $\mathbb{Z}^{+}$-semimodule, whence $Hom_{\mathbb{Z}^{+}}(S,D)$
is $N$-$i$-injective by Proposition \ref{dtohomd}.$\blacksquare $
\end{Beweis}

The following examples shows one of the advantages of Theorem \ref%
{c-i-injective}.

\begin{ex}
Let $S$ be an entire, cancellative, zerosumfree semiring. By Theorem \ref%
{c-i-injective}, every left $S$-semimodule $L$ can be embedded subtractively
in a $c$-$i$-injective left $S$-semimodule. On the other hand, if $L\neq 0,$
then $L$ cannot be embedded in an injective $S$-semimodule since the only
injective left $S$-semimodule is $\{0\}$ (by Lemma \ref{N0}). This is the
case, in particular, for $S:=\mathbb{Z}^{+}.$
\end{ex}


\end{document}